\documentclass[12pt]{article}

\usepackage{amsfonts,amsmath,bm,bbm}
\usepackage{thmtools}
\usepackage{thm-restate}
\usepackage{latexsym,color,epsfig,mathrsfs,enumerate}
\usepackage[toc,page]{appendix}
\usepackage{fancyhdr}
\usepackage[colorlinks=true, allcolors=red]{hyperref}
\usepackage{indentfirst}
\usepackage[maxbibnames=9]{biblatex}
\usepackage{comment}
\usepackage{xcolor}
\usepackage[shortlabels]{enumitem}

\usepackage{graphicx}
\usepackage{pgfplots}
\pgfplotsset{compat=1.15}
\usetikzlibrary{arrows}

\nocite{*}

\usepackage{amsthm}

\newtheorem{thm}{Theorem}[section]
\newtheorem{dfn}[thm]{Definition}
\newtheorem{lem}[thm]{Lemma}
\newtheorem{cor}[thm]{Corollary}
\newtheorem{prop}[thm]{Proposition}
\newtheorem{rmk}[thm]{Remark}

\newcommand{\mc}[1]{\mathcal{#1}}

\usepackage{combelow}

\usepackage{newunicodechar}
\newunicodechar{ș}{\cb{s}}
\newunicodechar{ț}{\cb{t}}

\newcommand{\eps}{\varepsilon}

\newcommand{\bP}{\mathbb P}
\newcommand{\bE}{\mathbb E}

\title{Isolation in Random Graphs}
\date{}

\author{Lauren\cb{t}iu Ploscaru\thanks{HUN-REN Alfr\'ed R\'enyi Institute of Mathematics, Budapest. Email: \texttt{laurentiuploscaru@renyi.hu},\\
supported by the ERC Advanced Grants no. 882971 GeoScape and ERMiD no. 101054936.
}}

\begin{document}

\maketitle

\begin{abstract}
We study isolation in the binomial random graph $G(n,p)$.
When $p=o(1)$ and $np\to\infty$, allowing an $F$-free graph
on the undominated vertices gives the same first-order
asymptotics for every fixed forest $F$ as the domination number $\gamma(G)$ in $G(n,p)$.
For the ordinary isolation number $\iota(G)$ ($F$ is an edge) we prove
two-point concentration when $p$ is bounded away from $1$ and $np^2/(\log n)^3\to\infty$. We also identify a transition in the gap $\gamma(G)-\iota(G)$ when $\log(1/p)=o(\log n)$. This gap tends to infinity in probability when $p$ is well below $(\log\log n/\log n)^2$ and
is at most one with high probability when $p$ is well above
this scale. Finally, for every fixed graph $F$, we show that
the classical minimum degree bound for domination is
asymptotically sharp for $F$-isolation, even for regular host graphs.
\end{abstract}

\section{Introduction}

The domination number is a central parameter in graph theory which has been studied extensively via an extremal perspective. A set $S\subseteq V(G)$ is \emph{dominating} in a graph $G$ if each of its vertices is either in $S$ or has a neighbor in $S$. The smallest size of such a set is called the \emph{domination number} of $G$, denoted by $\gamma(G)$. Equivalently, $\gamma(G)$ is the minimum number of closed neighborhoods needed to cover $V(G)$.\\ 
\indent Caro and Hansberg~\cite{CH} introduced isolation as a generalization of domination in which some vertices may be left undominated, provided that they induce no edges. For a vertex subset $S$ of a graph $G$, we let $N_G[S]$ denote its \emph{closed neighborhood}, i.e. $S$ together with all the vertices having a neighbor in $S$. A set $S\subseteq V(G)$ is \emph{isolating} if the subgraph induced on $U_S=V(G)\setminus N_G[S]$ is edgeless, and the minimum size of an isolating set is called the \emph{isolation number} $\iota(G)$ of $G$. This parameter is also known as the vertex-edge domination number as $N_G[S]$ touches every edge of $G$.\\
\indent More generally, let $\mc F$ be a nonempty family of graphs and $G$ a host graph. We say that a set $S\subseteq V(G)$ is \emph{$\mc F$-isolating} if $G-N_G[S]$ contains no member of $\mc F$ as a subgraph. We write
\[
 \iota(G,\mc F):=\min\{|S|:S\subseteq V(G)\text{ is $\mc F$-isolating}\}.
\]
\indent We will abbreviate $\iota(G,\{F\})$ to $\iota(G,F)$. In particular,
$\iota(G,K_1)=\gamma(G)$ and $\iota(G,K_2)=\iota(G)$. Since every dominating set is $\mc F$-isolating, we have $\iota(G,\mc F)\le\gamma(G)$.

The study of key parameters in random graphs is a central theme in probabilistic combinatorics, as it provides insight into their typical behavior for arbitrary graphs. Classical examples include the clique \cite{BE}, independence \cite{Frieze1990} and chromatic \cite{Bollobas1988} numbers, while the domination number has also been studied extensively~\cite{WG,GLS,BWZ}. In this paper, we initiate a systematic study of the isolation number in the binomial random graph model $G(n,p)$, where each edge is present independently with probability $p$. We take the vertex set of $G(n,p)$ to be $[n]=\{1,\ldots,n\}$. All logarithms are natural, and $f\ll g$ (or $g\gg f$)
means that $f/g\to0$. We use \emph{whp} to mean with high
probability, that is, with probability tending to one as
$n\to\infty$.

Our main aim is to understand how isolation compares with domination as the edge density varies.
We begin by showing that isolation and domination have the
same first-order asymptotics when $p=o(1)$ and $np\to\infty$. In fact, we prove the following more general result for
$F$-isolation, valid for every fixed forest $F$.

\begin{thm}\label{thm:sparse}
Let $F$ be a forest and let $p=p(n)$ satisfy $np\to \infty$ and $0<p\leq 1-n^{-o(1)}$. Then
$$\iota\left(G(n,p),F\right)
 =\frac{-\log(np)}{\log(1-p)}+O_F\left(1+\frac{\log\log(np)}{-\log(1-p)}\right)\quad whp.$$
In particular, if $p=o(1)$, then, whp
$\ \ \iota(G(n,p),F) \sim\gamma(G(n,p))
 \sim\dfrac{\log(np)}p$.
\end{thm}
For several classic random graph parameters, such as clique and independence
numbers at fixed edge density~\cite{BE}, and the chromatic
number in sufficiently sparse random graphs~\cite{AK}, a first-order estimate can be strengthened to \emph{two-point concentration}: with high probability, the parameter takes one of two consecutive integer values. The domination number of $G(n,p)$ follows the same behavior. Wieland and Godbole~\cite{WG} established it for fixed $p$, Glebov, Liebenau and Szab\'o~\cite{GLS} extended the range to $p\gg n^{-1/2}\log^2n$, and Bohman, Warnke and Zhu~\cite{BWZ} proved it for $p\ge n^{-2/3}\log^3n$. These results locate the domination number near the point where the expected number of dominating sets changes from tending to zero to tending to infinity.

We prove an analogous expectation criterion and two-point concentration result for isolation. Let $Y_k$ count the isolating $k$-sets in $G(n,p)$, and put $M_k=\bE [Y_k]$. Our result then says: 

\begin{thm}\label{thm:variable-counting}
Let $p=p(n)$ satisfy $np^2\log^{-3}(n)\to \infty$ and
$0<p\leq 1-\eta$ for some fixed $\eta>0$.
Let $k=k(n)$ be an integer satisfying $-k\log(1-p)=\log(np)-2\log\log(np)+O(1)$.
Then
\[
 \bP\bigl(\iota(G(n,p))\leq k\bigr)\longrightarrow
 \begin{cases}
  0, & \text{if } M_k\to0,\\
  1, & \text{if } M_k\to\infty.
 \end{cases}
\]
\indent Moreover, if $r$ is the smallest nonnegative integer for which
$M_r\geq 1/\log n$, then $r$ satisfies the same asymptotic
condition as $k$ above, and $\iota(G(n,p))\in\{r,r+1\}$ whp.
\end{thm}

The domination results above cover the whole range for $p$ in Theorem \ref{thm:variable-counting}. For isolation, the additional issue is to control the edges among the vertices missed by a candidate set. Our proof bounds the probability that two overlapping candidates are both isolating, which allows us to apply the second moment method at the isolation expectation threshold.

Theorem~\ref{thm:sparse} provides the same first-order estimate
for domination and isolation, but it does not determine
the size of their difference. Caro and
Hansberg~\cite{CH} asked how large
$\gamma(G)-\iota(G)$ can be among graphs of prescribed
order and minimum degree. Motivated by their question,
we next compare the two parameters more precisely in
$G(n,p)$. Our next result shows that the gap changes from tending to
infinity to being at most one at a specific threshold. 

\begin{thm}\label{thm:gap}
Let $\eta\in(0,1)$ and let $p=p(n)$ satisfy $-\log(p)=o(\log n)$ and
$0<p\leq1-\eta$. Setting $p_*:=(\log\log n/\log n)^2$, the following statements hold:
\begin{enumerate}[(i)]
\item If $p\ll p_*$, then whp $\gamma(G(n,p))-\iota(G(n,p))
 \sim 2p_*p^{-1}$.
\item If $p\gg p_*$, then whp $0\leq\gamma(G(n,p))-\iota(G(n,p))\leq1$.
\item For every fixed $c>0$, if $p\geq cp_*$, then whp $\gamma(G(n,p))-\iota(G(n,p))=O_{c,\eta}(1)$.
\end{enumerate}
\end{thm}

For example, taking $p=(\log n)^{-3}$ gives $\gamma(G(n,p))-\iota(G(n,p)) 
\sim 2\log n\,(\log\log n)^2\ $ whp,
even though the ratio of the two parameters tends to one.
These estimates describe the gap in typical random graphs,
whereas the question of Caro and Hansberg concerns its
largest possible value for given order and minimum degree.
At the transition $p\sim cp_*$, a finer question is to
determine which integer values the difference takes
and with what probabilities.

Finally, we apply our random graph estimates to bound the isolation number in terms of the minimum
degree. For graphs $G$ with $\delta(G)\geq d$, the classical domination bound gives
\[
 \iota(G,F)\leq\gamma(G)\leq\frac{\log(d+1)+1}{d+1}\cdot |V(G)|.
\]
\indent In the case of ordinary isolation, Li, Zhang and Ye~\cite{LZY} obtained a smaller additive constant in the numerator. Caro and Hansberg~\cite{CH} gave a construction showing that the leading coefficient in front of $(\log d)/d$ must be at least $1/2$.
Borg and Caro~\cite{BC} extended this lower bound to every fixed graph $F$ by showing that, for every
fixed graph $F$ and every sufficiently large integer $d$, there exists a graph $G$ with $\delta(G)\geq d$ and $\iota(G,F)\geq (1-o(1))\log d\cdot |V(G)| \bigl(d|V(F)|\bigr)^{-1}$. They also obtained this coefficient
using regular host graphs when $F$ is regular.

Therefore, the general lower and upper bounds differed
by a factor of $|V(F)|$ in their leading terms.
We close this gap by showing that the coefficient $1$
is best possible for every fixed $F$, even when
the host graph is required to be regular.

\begin{thm}\label{thm:extremal}
Let $F$ be a fixed non trivial graph.
As $d\to\infty$ through the integers, we have\vspace{3pt}
$\phantom{x}\qquad\displaystyle\sup_{\delta(G)\geq d}
 \dfrac{\iota(G,F)}{|V(G)|}
 =\dfrac{\log d+O_F(\log\log d)}{d}$.
Moreover, $\displaystyle\sup_{G\text{ is }d\text{-regular}}
 \dfrac{\iota(G,F)}{|V(G)|}
 \sim\dfrac{\log d}{d}$.
\end{thm}
Theorem~\ref{thm:sparse}
gives a sharper bound for forests $F$: there are graphs $G$
with $\delta(G)\geq d$ and
\[
 \frac{\iota(G,F)}{|V(G)|}
 \geq\frac{\log d-2\log\log d-O_F(1)}{d}.
\]
\indent To obtain the regular construction, we couple a random
$d$-regular graph between two binomial random graphs.
The larger graph controls the neighborhoods of candidate
isolating sets, while the smaller graph supplies copies
of $F$ outside those neighborhoods.

The remaining of the paper is structured as follows. Section~\ref{sec:preliminaries} collects the estimates used in the proofs. In Section~\ref{sec:first-order}, we prove the
first-order estimate for forests by showing that small sets
leave enough undominated vertices to contain a copy of $F$.
We use the first and second moments of the number of isolating
sets to prove two-point concentration in
Section~\ref{sec:concentration}, then compare the expectation
thresholds for domination and isolation in
Section~\ref{sec:comparison}. Next, in Section~\ref{sec:regular}, we
give binomial and regular constructions showing that the leading
coefficient in the minimum degree bound is sharp.
We finish with two questions in Section~\ref{sec:concluding-remarks}.

\section{Preliminary Tools}\label{sec:preliminaries}

We begin with some numerical estimates for some functions, which we state below.

\begin{lem}
    For all integers $n$ and $k$ with $1\leq k\leq n$, we have
$\left(\dfrac{n}{k}\right)^k
 \leq \dbinom{n}{k}
 \leq \left(\dfrac{en}{k}\right)^k$.
\end{lem}

\begin{lem}\label{log-1-p-approx}
    Let $p\in (0,1)$. Then $p\leq -\log(1-p)\leq p-p\log(1-p)$.
\end{lem}

\begin{proof}
    Setting $q:=-\log(1-p)$, the first inequality becomes $1-e^{-q}\leq q$, which is Bernoulli. Conversely, the second one rewrites as $q\leq (1+q)(1-e^{-q})$, which reduces to $(1+q)e^{-q}\leq 1$, or $e^q\geq 1+q$.
\end{proof}

\begin{lem}\label{lower-product-approx}
    Let $u_1,u_2,\ldots,u_n\in[0,1]$, where $n\geq 1$ is an integer. Then
$\displaystyle\prod_{i=1}^n(1-u_i)\geq1-\displaystyle\sum_{i=1}^n u_i$.
\end{lem}

\begin{proof}
    We proceed by induction on $n$. The case when $n=1$ is trivially true, so fix $n$ and let as assume the conclusion holds true for any $n$ numbers in $[0,1]$. Let $u_1,u_2,\ldots,u_{n+1}\in [0,1]$ and write $S:=\sum_{i=2}^{n+1}u_i$. By induction, we have $\prod_{i=1}^{n+1}(1-u_i)\geq (1-u_1)(1-S)=1-u_1-S+Su_1\geq 1-\sum_{i=1}^{n+1}u_i$, which is what we needed to prove.   
\end{proof}

We will also require some classic concentration inequalities. See e.g. appendix A in \cite{alon04}.

\begin{thm}[\textbf{Chernoff Inequality}]\label{cher}
Let $X$ be a random variable with binomial distribution and let $\mu = \bE[X]$. Then, for $0\leq \delta\leq 1$, the following inequalities hold:
$$ \bP\big(X\leq (1-\delta )\mu \big)\leq \exp\left(-{\dfrac {\delta ^{2}\mu }{2}}\right).$$
$$ \bP\big(X\geq (1+\delta )\mu \big)\leq \exp\left(-{\dfrac {\delta ^{2}\mu }{4}}\right).$$
\end{thm}

\begin{thm}[\textbf{Markov Inequality}]\label{markov}
Let $X$ be a nonnegative random variable with finite expectation
$\mu=\bE[X]$. Then $\bP\big(X\geq t\big)\leq \mu t^{-1}$ for every $t>0$.
\end{thm}

Finally, we will make use of the following binomial lower-tail bound. 

\begin{lem}\label{lem:binomial-lower-tail}
Let $X\sim\operatorname{Bin}(N,p)$, where $N\geq1$ is an
integer and $0<p<1$. Then, for every integer $L$ satisfying
$1\leq L\leq Np$, the following holds:
\begin{align*}
 \bP(X\leq L)
 &\leq (L+1)(1-p)^N
       \left(\frac{eNp}{L(1-p)}\right)^L.
\end{align*}
\end{lem}

\begin{proof}
We start by comparing the probabilities that $X$ takes consecutive values. Notice that for $0\leq j<L$, we have:
\[
 \frac{\bP(X=j+1)}{\bP(X=j)}
 =\frac{(N-j)p}{(j+1)(1-p)}
 \geq 1.
\]
\indent Consequently, when we write $\bP(X\leq L)$ as a sum, we can upper bound each $\bP(X=j)$ with $0\leq j\leq L$ by $\bP(X=L)$. We therefore get:  
$$\bP(X\leq L)\leq (L+1) \binom NL\cdot p^L(1-p)^{N-L}
 \leq (L+1)(1-p)^N      \left(\frac{eN}{L}\cdot \frac{p}{(1-p)}\right)^L.$$
\end{proof}

\section{First order estimate of the isolation number}\label{sec:first-order}

We prove Theorem~\ref{thm:sparse} in two steps. First, we show that a large enough vertex set $S$ leaves very few undominated vertices. For the lower bound, we compare the number of vertices left undominated by a set $S$ with the size of an $F$-free set. The two lemmas below make this comparison explicit: if $S$ is too small, there will be enough vertices left to span $F$.

\begin{lem}\label{lem:expansion}
Let $G\sim G(n,p)$, where $p=p(n)$ satisfies $np\to \infty$ and $0<p\leq 1-n^{-o(1)}$. 
For each $S\subseteq V(G)$, let $U_S$ denote the set of
vertices outside $S$ having no neighbor inside $S$.
Then, for every $\eps=\eps(n)\in(0,1)$ satisfying $(np)^\eps\geq e^4\log^2(np)$, whp every set $S\subseteq V(G)$ with
$|S|\le (1-\eps)\log(np)/(-\log(1-p))$
also satisfies $|U_S|\ge (np)^\eps/(3p)$.
\end{lem}

\begin{proof}
Put $k=\lfloor(1-\eps)\log(np)/(-\log(1-p))\rfloor$ and fix a set $S$ of size $k$. Since $p\leq -\log(1-p)$ by Lemma \ref{log-1-p-approx}, we have $k=o(n)$. 
We note that $|U_S|\sim\operatorname{Bin}\left(n-k,(1-p)^k\right)$. Indeed, each vertex outside $S$ belongs to $U_S$ precisely when all its $k$ edges to $S$ are absent. This event has probability $(1 - p)^k$, independently for different vertices.

It follows that $\bE|U_S|=(n-k)(1-p)^k
 \ge(np)^\eps/(2p)$ for all large enough $n$. Then by Chernoff's inequality with $\delta=1/3$ we get
\[
 \bP\left(|U_S|<\frac{(np)^\eps}{3p}\right)
 \le\exp\left(-\frac{(np)^\eps}{36p}\right).
\]
\indent We now wish to take a union bound over all $\binom{n}{k}$ choices of $S$. 
Since $k\leq \log(np)/p=o(n)$ and the function
$x\mapsto x\log(en/x)$ is increasing on $(0,n)$, we have, for large enough $n$:
\[
 \log\binom nk
 \leq k\log\left(\frac{en}{k}\right)
 \leq \frac{\log(np)}p
       \log\left(\frac{enp}{\log(np)}\right)
 \leq \frac{\log^2(np)}p.
\]
\indent  The probability that there is a set $S$ of size $k$ with $|U_S|<\dfrac{(np)^\eps}{3p}$ is then at most
\[
 \binom nk\exp\left(-\frac{(np)^\eps}{36p}\right)
 \leq
 \exp\left(\frac{\log^2(np)}p-\frac{(np)^\eps}{36p}\right)
 =o(1).
\]
\indent The hypothesis $(np)^\eps\geq e^4\log^2(np)$ is precisely what makes the last bound hold true.\\
\indent Finally, any set $T$ of size smaller than $k$ can be extended
to a set $S$ of size $k$. Since $U_S\subseteq U_T$, we also have
$|U_T|\geq |U_S|\geq (np)^\eps/3p$,
completing the proof.
\end{proof}

\begin{lem}\label{lem:forest-free}
Let $G\sim G(n,p)$, where $p=p(n)$ satisfies
$np\to\infty$ and $0<p\leq 1-n^{-o(1)}$.
Let $F$ be a fixed forest with at least one vertex.
Then whp every set of
$\left\lceil 3\log(np)/(-\log(1-p))\right\rceil$
vertices in $G$ spans a copy of $F$.
\end{lem}

\begin{proof}
The case when $F$ has only one vertex is immediate. So assume $r:=|V(F)|\ge2$ and set $t=\left\lceil 3\log(np)/(-\log(1-p))\right\rceil$. Any graph of minimum degree at least $r-1$ contains $F$, by greedily embedding its tree components. Thus every nonempty subgraph of an $F$-free graph must have a vertex of degree at most $r-2$. Successively deleting such vertices shows that an $F$-free graph on $t$ vertices has at most $(r-2)t< rt$ edges.

We show that whp no set of $t$ vertices spans fewer than $rt$
edges. Write $q:=-\log(1-p)$ for convenience. Our assumption on $p$ gives $q=o(\log(np))$: if $p\leq 1/2$ then $q\leq \log 2$, while if $p > 1/2$ then $\log(np)\sim \log n$. Besides, note that $qt=(3+o(1))\log(np)$, while by Lemma~\ref{log-1-p-approx},
\begin{equation}\label{eqn:p-q-bounds}
     3\log(np)+q\ge qt\ge pt\geq \frac{3p\log(np)}q
 \geq \frac{3\log(np)}{1+q}.
\end{equation}
\indent In particular $pt\to \infty$ as $n\to \infty$, and so $rt< p\binom t2$ for all large enough $n$. For a fixed set $T$ of size $t$, the number of edges in $G[T]$ has distribution $\operatorname{Bin}\bigl(t(t-1)/2,p\bigr)$. Applying Lemma~\ref{lem:binomial-lower-tail} with
$N=\binom t2$ and $L=rt-1$, we obtain
\begin{align*}
 \bP\bigl(G[T]\text{ is }F\text{-free}\bigr)
 &\leq \bP\bigl(e(G[T])<rt\bigr)\leq rt(1-p)^{\binom t2}
       \left(\frac{ep t(t-1)/2}
       {(rt-1)(1-p)}\right)^{rt-1}\\
 &\leq rt\,e^{-q\binom t2}
       \left(\frac{ept}{1-p}\right)^{rt}=\exp\left(-q\binom t2+\log(rt)
       +rt\bigl(1+\log(pt)+q\bigr)\right).
\end{align*}
\indent Here we used $\binom t2/(rt-1)\leq t$ and the fact that
$ept/(1-p)>1$ for all large enough $n$.\\
\indent We next estimate the factors accompanying $e^{-q\binom t2}$. From (\ref{eqn:p-q-bounds}) we get $\log(pt)=O(\log\log(np))$. As $r$ is fixed and
$q=o(\log(np))$, it follows that $\log(rt)+rt\bigl(1+\log(pt)+q\bigr)
 =o\bigl(t\log(np)\bigr)$. Moreover, 
$qt(t-1)/2 =\left(3/2+o(1)\right)\cdot t\log(np)$. It follows that
\[
 \bP\bigl(G[T]\text{ is }F\text{-free}\bigr)
 \leq \exp\mbox{\large(}-(3/2+o(1))\cdot t\log(np)\mbox{\large)}.
\]
\indent Finally, since $pt\to\infty$, for all large enough $n$,
 $\dbinom nt\leq \left(\dfrac{en}{t}\right)^t \leq (np)^t$.\\
\indent The expected number of $F$-free sets of size $t$ is therefore
at most
\[
 (np)^t\exp\left(-(3/2+o(1))t\log(np)\right)
 =\exp\left(-(1/2+o(1))t\log(np)\right)
 =o(1).
\]
\indent The conclusion follows from Markov's inequality, i.e. the first moment method.
\end{proof}

\begin{proof}[Proof of Theorem~\ref{thm:sparse}]
Write $q:=-\log(1-p)$ for convenience. As in the proof of
Lemma~\ref{lem:forest-free}, we have $q=o(\log(np))$.
For the upper bound, fix a set $S$ of size
$\left\lceil\ q^{-1}\left(\log(np)+\log\log(np)\right)\right\rceil$,
and let $U_S$ be the set of vertices outside $S$ having
no neighbor in $S$. As in the
proof of Lemma~\ref{lem:expansion}, we have
 $\bE|U_S|=(n-|S|)(1-p)^{|S|}
 \leq n e^{-q|S|}
 \leq (p\log(np))^{-1}$.
By Markov's inequality:
\[
  \bP\left(|U_S|>\frac1q\right)
 \leq \frac{q}{p}\cdot\log^{-1}(np)
 \leq \frac{1+q}{\log(np)}
 =o(1).
\]
\indent Adding $U_S$ to $S$ gives a dominating set, which is also
$F$-isolating. It follows whp that
\[
 \iota(G,F)\leq\gamma(G)\leq |S|+|U_S|
 \leq \frac{\log(np)+\log\log(np)+1+q}{q}.
\]
\indent Turning to the lower bound, put $\eps=(2\log\log(np)+4)/\log(np)$.
For all large enough $n$, we have $0<\eps<1$ and
$(np)^\eps=e^4\log^2(np)$.
Lemma~\ref{lem:expansion} shows that whp every set $S$ whose size is
$|S|\leq q^{-1}(\log(np)-2\log\log(np)-4)$ satisfies
\[
 |U_S|\geq\frac{e^4\log^2(np)}{3p}
 >\frac{3\log(np)}p+1
 \geq\left\lceil\frac{3\log(np)}q\right\rceil,
\]
where the last inequality uses
Lemma~\ref{log-1-p-approx}.
Combining this with Lemma~\ref{lem:forest-free}, we conclude
that whp $G[U_S]$ contains a copy of $F$ for every such $S$.
Thus none of these sets is $F$-isolating, so
\[
 \iota(G,F)>\frac{\log(np)-2\log\log(np)-4}{q}\quad whp.
\]
\indent Together with the upper bound, this proves the first
assertion.\\
\indent Finally, when $p=o(1)$, Lemma~\ref{log-1-p-approx} gives
$q\sim p$. Taking $F=K_1$ (a vertex) and $F=K_2$ (an edge), and recalling that
$\iota(G,K_1)=\gamma(G)$ and $\iota(G,K_2)=\iota(G)$, gives the last part as well, completing the proof.
\end{proof}

We remark that the forest assumption is essential in Theorem \ref{thm:sparse}. Take $F=K_3$ and set $p=n^{-4/5}$ for example. Then the expected number of triangles is at most $(pn)^3=O(n^{3/5})$. Taking one vertex from each triangle gives a $K_3$-isolating set, therefore by Markov's inequality we get whp that $\iota(G(n,p),K_3)\le n^{7/10}=o(\log(np)/p)$.
Thus a cyclic forbidden graph can have a much smaller isolation number in the sparse range.

\section{Two-point concentration}\label{sec:concentration}

Recall that $Y_k$ counts the isolating sets of size $k$ in $G(n,p)$ and that $M_k=\bE[Y_k]$. We prove Theorem~\ref{thm:variable-counting} by estimating the first and second moments of $Y_k$. Markov's inequality already gives the nonexistence of an isolating set of size at most $k$ when $M_k\to0$, so the main step is to show that whp there is an isolating set of size $k$ when $M_k\to\infty$. We then obtain the two-point conclusion by comparing sets of consecutive sizes.

Throughout this section, let $\eta\in(0,1)$ be fixed and whenever we refer to $G(n,p)$, we assume that $p(n)$ satisfies the conditions 
\begin{equation}\label{eq:s4-range}
 0<p=p(n)\leq 1-\eta,
 \qquad np^2\log^{-3}(n)\longrightarrow\infty.
\end{equation}
\indent For convenience, we will also set $q:=-\log(1-p)$ and $ L:=\log(np)$ throughout this section.

Let us note that, by symmetry, if $S$ is any fixed set of $k$ vertices, then $M_k=\binom{n}{k}\cdot \bP(S\text{ isolates})$. 
We begin by estimating this probability for $1\leq k<n$. As in Lemma~\ref{lem:expansion}, let $U_S$ denote the set of vertices outside $S$ having no neighbor in $S$. Write $a:=(1-p)^k$, so that $|U_S|\sim\operatorname{Bin}(n-k,a)$. To determine $U_S$, we expose only edges between $S$ and its complement. Conditional on $U_S$, the edges within $U_S$ remain independent, each present with probability $p$.

The set $S$ is isolating precisely when $U_S$ is independent. Conditional on $|U_S|=t$, this happens with probability
$(1-p)^{\binom t2}$. Summing over all possible sizes of $U_S$ gives
\begin{align}\label{eq:isolation-exact-probability}
 \bP(S\text{ isolates})
 &=\sum_{t=0}^{n-k}\binom{n-k}{t}a^t(1-a)^{n-k-t}
       (1-p)^{\binom t2}\notag\\
 &=(1-a)^{n-k}\sum_{t=0}^{n-k}\binom{n-k}{t}
       \left(\frac{a}{1-a}\right)^t(1-p)^{\binom t2}.
\end{align}
\indent To compare isolation with domination, we need to account for the independent
set of vertices that may remain undominated. We remark that the factor $(1-a)^{n-k}$ is exactly the probability that $S$ is dominating, which also corresponds to the term $1$ (for $t=0$) in the sum; the other terms account for the extra possibilities allowed by isolation. We also note that for values of $t$ that contribute significantly to (\ref{eq:isolation-exact-probability}) we can make the following approximation: $$\binom{n-k}{t}\left(\frac{a}{1-a}\right)^t
 \approx\frac{1}{t!}\left(\frac{(n-k)a}{1-a}\right)^t.$$
\indent It is therefore natural to introduce the following potential potential, which will aid us in quantifying the resulting gain of isolation in the first-moment calculation.

\begin{dfn}\label{def:isolation-potential}
For $\lambda>0$, we define the \textbf{isolation} $\lambda$-\textbf{potential} by
\[
 \Phi_\lambda(x):=
 \sum_{t\geq0}\frac{x^t}{t!}
 e^{-\lambda\binom t2}
 \qquad(x\geq0).
\]
\end{dfn}
In the applications below, we will take $\lambda=q$, so that
$e^{-\lambda\binom t2}=(1-p)^{\binom t2}$.\\ 
\indent We first study $\Phi_\lambda$ for a general $\lambda>0$.
The ratio of consecutive summands is $xe^{-\lambda t}/(t+1)$.
When $\lambda$ stays bounded and $\lambda x\to\infty$, this
suggests that the main contribution comes from values of $t$ of
order $\log(\lambda x)/\lambda$. The next lemma gives the growth
and tail estimates that we will use to approximate the finite sum
in~\eqref{eq:isolation-exact-probability}. These estimates will
also be useful when bounding the probability that two fixed sets
are both isolating.

\begin{lem}\label{lem:series}
Fix $\lambda_0>0$. Let $0<\lambda\leq\lambda_0$ and suppose
that $\lambda x\to\infty$. Then
\begin{equation*}
 \log\Phi_\lambda(x)
 =\frac{\log^2(\lambda x)
 +O_{\lambda_0}\bigl(\log(\lambda x)\log\log(\lambda x)\bigr)}
 {2\lambda}
 \qquad\text{and}\qquad
 \frac{x\cdot \Phi_\lambda'(x)}{\Phi_\lambda(x)}
 =O_{\lambda_0}\left(\frac{\log(\lambda x)}{\lambda}\right).
\end{equation*}
\indent Moreover, for all sufficiently large $\lambda x$, the following estimates uniformly hold:
\begin{equation*}
 \sum_{t=0}^{\lfloor\frac{1}{2\lambda}\log(\lambda x)\rfloor}
 \frac{x^t}{t!}e^{-\lambda\binom t2}
 \leq\Phi_\lambda(x)^{5/6},
\qquad\text{and}\qquad\label{eq:s4-series-upper-tail}
 \sum_{\lambda t>16\log(\lambda x)}
 (1+t)\frac{x^t}{t!}e^{-\lambda\binom t2}
 =o(1).
\end{equation*}
\end{lem}

\begin{proof}
Put $R:=\log(\lambda x)$ and, for each integer $t\geq0$, write
$w_t(x):=x^t\cdot e^{-\lambda\binom t2}/(t!)$.
Using the elementary bounds $(t/e)^t\leq t!\leq t^t$
for $t\geq1$, we deduce that
\begin{equation}\label{eq:s4-summand}
 \log w_t(x) \leq t\left(R+1+\frac{\lambda}{2}
                  -\log(\lambda t)\right)
       -\frac{\lambda t^2}{2}
 \leq t\left(R+1+\frac{\lambda}{2}\right)
       -\frac{\lambda t^2}{2}+\frac{1}{e\lambda},
\end{equation}
where the last inequality follows from $-u\log u\leq1/e$
for $u>0$.

We first bound the upper tail. If $\lambda t>16R$, then,
for all large enough $R$, we have $\log(\lambda t)\geq0$
and $R+1+\lambda/2\leq2R$, since $\lambda\leq\lambda_0$.
Therefore, the first inequality in~\eqref{eq:s4-summand} gives us $\log w_t(x)\leq$
 $t\left(4R-\lambda t\right)/2 \leq-\lambda t^2/4 \leq-t$.
Consequently, $\displaystyle\sum_{\lambda t>16R}(1+t)w_t(x)
 \leq\displaystyle\sum_{t>16R/\lambda_0}(1+t)e^{-t}
 =o(1)$. This proves the upper-tail estimate, uniformly for $0<\lambda\leq\lambda_0$.

We next estimate the full sum. Completing the square in
the last line of~\eqref{eq:s4-summand} gives
\[
 \log w_t(x)
 \leq\frac{(R+1+\lambda/2)^2}{2\lambda}
       +\frac{1}{e\lambda}
 =\frac{R^2}{2\lambda}
       +O_{\lambda_0}\left(\frac{R}{\lambda}\right).
\]
\indent There are at most $1+16R/\lambda$ terms with
$\lambda t\leq16R$. Combining this bound with the
upper-tail estimate, and using
$\log(1+16R/\lambda)\leq16R/\lambda$, we obtain that 
 $\log\Phi_\lambda(x)\leq R^2/(2\lambda)+O_{\lambda_0}\left(R/\lambda\right)$.

For the reverse inequality, it is enough to retain one
summand. Put $m:=\lfloor R/\lambda\rfloor$, which is positive
for all large enough $R$. Using $\lambda m\leq R$, we get $\log\Phi_\lambda(x)\geq\log w_m(x)$ and so
\begin{equation*}
 \log\Phi_\lambda(x)\geq mR-\frac{\lambda m^2}{2}-m\log R
 =\frac{R^2}{2\lambda}-\frac{\lambda}{2}\left(m-\frac{R}{\lambda}\right)^2-m\log R
 \geq\frac{R^2}{2\lambda}-\frac{R\log R}{\lambda}-\frac{\lambda}{2}.
\end{equation*}
\indent Together with the upper bound, this proves the claimed
estimate for $\log\Phi_\lambda(x)$.

We now turn to the lower tail. The quadratic in the rightmost side  of~\eqref{eq:s4-summand} is increasing for
$0\leq t\leq R/(2\lambda)$. Hence, throughout this range, we have
$\log w_t(x) \leq 3R^2/(8\lambda) +O_{\lambda_0}\left(R/\lambda\right)$.
Summing over at most $1+R/(2\lambda)$ terms and comparing
with our estimate for $\log\Phi_\lambda(x)$ yields
\[
 \log\left(\sum_{0\leq\lambda t\leq R/2}w_t(x)\right)
 \leq\frac{3R^2}{8\lambda}
       +O_{\lambda_0}\left(\frac{R}{\lambda}\right)
 =\left(\frac34+o(1)\right)\log\Phi_\lambda(x).
\]
\indent Since $3/4<5/6$, exponentiating gives
$\displaystyle\sum_{0\leq\lambda t\leq R/2}w_t(x)\leq\Phi_\lambda(x)^{5/6}$
for all large enough $R$.

Finally, $\Phi_\lambda(x)$ can be differentiated term
by term, giving
 $\dfrac{x\cdot \Phi_\lambda'(x)}{\Phi_\lambda(x)}
 =\dfrac{\sum_{t\geq0}t\cdot w_t(x)}{\Phi_\lambda(x)}$.

Since $\Phi_\lambda(x)\geq1$, splitting the numerator at $t=16R/\lambda$ and applying the upper-tail estimate, we obtain
\[
 \frac{x\cdot \Phi_\lambda'(x)}{\Phi_\lambda(x)}
 \leq\frac{16R}{\lambda}
      +\frac{1}{\Phi_\lambda(x)}\cdot \displaystyle\sum_{\lambda t>16R}t\cdot w_t(x)
 \leq\frac{16R}{\lambda}+o(1)
 =O_{\lambda_0}\left(\frac{R}{\lambda}\right).
\]
\indent All estimates above hold uniformly for
$0<\lambda\leq\lambda_0$, completing the proof.
\end{proof}

We apply Lemma~\ref{lem:series} with $\lambda=q=-\log(1-p)$.
Since $p\leq1-\eta$, we have $0<q\leq-\log\eta$, so its
estimates apply with $\lambda_0=-\log\eta$.
The next lemma justifies replacing the finite sum
in~\eqref{eq:isolation-exact-probability} by the isolation
potential, comparing the probabilities that a
set is isolating or dominating.

\begin{lem}\label{lem:poissonization}
Fix $\eta\in(0,1)$ and let $G\sim G(n,p)$, where
$0<p\leq1-\eta$ and $np^2/\log^{3}n\to\infty$.
Let $S\subseteq[n]$ be a fixed set of size $k=k(n)$ satisfying
$kq=\log(np)-2\log\log(np)+O(1)$. Then:
\begin{align*}
 \bP(S\text{ dominates})
 &=\bigl(1-(1-p)^k\bigr)^{n-k},\\
 \bP(S\text{ isolates})
 &=\bigl(1-(1-p)^k\bigr)^{n-k}\cdot
   \Phi_q\left((n-k)(1-p)^k\right)(1+o(1)),\\
 \log\Phi_q\left((n-k)(1-p)^k\right)
 &=(2+o(1))\cdot\frac{(\log\log(np))^2}{q}.
\end{align*}
\indent The $o(1)$ errors are uniform in $S$ and in $k$ whenever the
$O(1)$ term in the hypothesis ranges over a fixed bounded
interval.
\end{lem}

\begin{proof}
Set $L:=\log(np)$, $a:=(1-p)^k$ and $\mu:=a(n-k)$. Observe that $p\leq q\leq(1-\log\eta)p$ by Lemma~\ref{log-1-p-approx}.
Together with the hypothesis, we get
$k=\Theta_\eta(L/p)=o(n)$, $a=\Theta\left(L^2/(np)\right)=o(1)$ and 
$q\mu=\Theta_\eta(L^2)$, which implies that $q\mu\to\infty$, so Lemma~\ref{lem:series} does apply.

Recall that $|U_S|\sim\operatorname{Bin}(n-k,a)$, so $\mu=\bE[|U_S|]$.
Since $S$ is dominating precisely when $U_S=\emptyset$, the first equality we need to prove follows immediately.

For isolation, we go back to~\eqref{eq:isolation-exact-probability}.
Put $x:=\mu/(1-a)$ and note that
$\dbinom{n-k}{t}\left(\dfrac{a}{1-a}\right)^t \leq\dfrac{x^t}{t!}$.

It follows that $\bP(S\text{ isolates})\leq (1-a)^{n-k}\cdot \Phi_q(x)$.
We show that the reverse inequality holds up to a relative
error $o(1)$ by restricting to the terms that contribute
most of the sum.

Set $T:=\lfloor16\log(qx)/q\rfloor$. Uniformly, for any positive integer $t\leq T$ the following holds: 
\[
 \binom{n-k}{t}\left(\frac{a}{1-a}\right)^t
 =\frac{x^t}{t!}\cdot
   \prod_{i=1}^{t-1}\left(1-\frac{i}{n-k}\right)
 =\bigl(1-o(1)\bigr)\frac{x^t}{t!}.
\]
To prove the last equality, we only need to lower bound the product. Indeed, by Lemma \ref{lower-product-approx}, 
$$\prod_{i=1}^{t-1}\left(1-\frac{i}{n-k}\right)\geq 1-\dfrac{\sum_{i=1}^{t-1}i}{n-k}>1-\dfrac{T^2}{2(n-k)}=1- O_\eta\left(\frac{(\log L)^2}{np^2}\right)
 =1-o(1).$$
\indent Moreover, the upper-tail estimate in Lemma~\ref{lem:series}
gives $\displaystyle\sum_{t>T}\dfrac{x^t}{t!}e^{-q\binom t2}=o(1)$.
As $\Phi_q(x)\geq1$, summing the term-wise approximation shows that the sum in~\eqref{eq:isolation-exact-probability} is $(1+o(1))\cdot \Phi_q(x)$. Thus
\[
 \bP(S\text{ isolates})
 =(1+o(1))_{\ }(1-a)^{n-k}\cdot \Phi_q(x).
\]
\indent It remains to replace $x$ by $\mu$ in the equality above. For every
$u\in[\mu,x]$, we have $qu=\Theta_\eta(L^2)$.
Integrating the derivative estimate in
Lemma~\ref{lem:series}, we obtain
\[
 0\leq\log\Phi_q(x)-\log\Phi_q(\mu)
 =\int_\mu^x\frac{\Phi_q'(u)}{\Phi_q(u)}\,du
 =O_\eta\left(\frac{\log L}{q}\log\frac{x}{\mu}\right).
\]
\indent Since $\log(x/\mu)=-\log(1-a)=O(a)$, the last expression
is $ O_\eta(a\log L/q)=O_\eta\bigl(L^2\log L/(np^2)\bigr)$, which is $o(1)$.
 It follows that $\Phi_q(x)=(1+o(1))\cdot \Phi_q(\mu)$,
proving isolation equality.

Finally, $q\mu=\Theta_\eta(L^2)$ implies
$\log(q\mu)=2\log L+O_\eta(1)$. Lemma~\ref{lem:series} now gives
\[
 \log\Phi_q(\mu)
 =\frac{\log^2(q\mu)
   +O_\eta\bigl(\log(q\mu)\log\log(q\mu)\bigr)}{2q}
 =(2+o(1))\frac{(\log L)^2}{q},
\]
which proves the last claim, thus completing the proof.
\end{proof}

Recall $M_k=\binom nk\cdot\bP(S\text{ isolates})$. By Lemma~\ref{lem:poissonization} we get, when
$qk=L-2\log L+O(1)$, that
\begin{equation*}
 \log M_k=\log\binom nk
 +(n-k)\log\bigl(1-(1-p)^k\bigr)
 +\log\Phi_q\bigl((n-k)(1-p)^k\bigr)+o(1).
\end{equation*}
\indent The first two terms give the logarithm of the expected number of
dominating $k$-sets. The additional term records the gain from
allowing an independent set of undominated vertices.

To apply the second moment method, we need to estimate the probability
that two fixed sets are both isolating. Sets with a small intersection
will behave almost independently. Sets that differ in few vertices
are more strongly dependent, but there are far fewer such pairs.
The next lemma gives the bounds we will use for these different
possibilities.

\begin{lem}\label{lem:pair}
Fix $\eta\in(0,1)$ and let $G\sim G(n,p)$, where
$0<p\leq1-\eta$ and $np^2/\log^3 n\to\infty$.
Let $S,T\subseteq[n]$ be distinct fixed sets of size $k=k(n)$,
where $kq=\log(np)-2\log\log(np)+O(1)$. Setting $\mu:=(n-k)(1-p)^k$,
$\rho:=\bP(S\text{ isolates})$ and $\ell:=|S\cap T|$, then: 
\begin{align*}
 \bP(S,T\text{ both isolate})
 &\leq\rho^2\exp\bigl(4\mu e^{-q(k-\ell)}+o(1)\bigr),\\
 \bP(S,T\text{ both isolate})
 &\leq\rho_{\ } \exp\left(-\bigl(1-o(1)\bigr)\cdot \mu(1-e^{-q(k-\ell)})\right).
\end{align*}
\indent Moreover, if $q\ell\leq1$, we have $\bP(S,T\text{ both isolate})
 \leq\rho^2\exp\bigl(o(\ell+1)\bigr)$.
 
 For any fixed bound on the $O(1)$ term, all the $o(1)$ error terms are
uniform in $k,S,T$.
\end{lem}

\begin{proof}
Once again, set $L:=\log(np)$ and $a:=(1-p)^k$. We have
$p\leq q\leq(1-\log\eta)p$, as in Lemma~\ref{lem:poissonization}, which leads to $k=\Theta_\eta(L/p)=o(n)$, $a=\Theta\left(L^2/(np)\right)=o(1)$ and 
$q\mu=\Theta_\eta(L^2)$. In particular, we have $ka=O_\eta(L^3/(np^2))=o(1)$.

Consider a vertex $v$ outside $S\cup T$ which is undominated
by $S$. For $v$ to be undominated by $T$ as well, its
$k-\ell$ edges to $T\setminus S$ must also be absent.
We therefore put
\[
 \theta:=\bP(v\in U_T\mid v\in U_S)
        =(1-p)^{k-\ell}.
\]
\indent Thus $v$ is undominated by both sets with probability
$a\theta$, and undominated only by $S$, or only by $T$,
with probability $a(1-\theta)$ in either case.
The remaining possibility is that $v$ is dominated by
both sets, which has probability $b:=1-a(2-\theta)$.
These choices are independent for different vertices
outside $S\cup T$, and depend only on their edges to
$S\cup T$.

Suppose that the common undominated group has size $j$,
and the two groups undominated by only one set have sizes
$u$ and $v$, respectively. For both sets to isolate, all
$\binom j2+\binom u2+\binom v2+j(u+v)$
possible edges within the two undominated sets must be
absent. There is no restriction on edges between the two
exclusive groups. 

There are $n-2k+\ell$ vertices outside $S\cup T$.
Write $x:=(n-2k+\ell)\cdot ab^{-1}$, expand over the three group
sizes and use
$(n-2k+\ell)! \leq(n-2k+\ell)^{j+u+v} (n-2k+\ell-j-u-v)!$ to get
\begin{align*}\label{eq:s4-joint-partition}
 \bP(S,T\text{ both isolate})
 \leq b^{n-2k+\ell}\cdot J(x,\theta),\quad \text{where}\\ 
 J(x,\theta):=
 \sum_{i\geq0}\frac{(x\theta)^i}{i!}e^{-q\binom i2}\cdot
 \Phi_q\bigl(x(1-\theta)e^{-qi}\bigr)^2.
\end{align*}
\indent To obtain the two copies of $\Phi_q$, fix the common
group of size $j$. Each vertex in an exclusive group must
avoid its $j$ vertices, giving the factor $e^{-qj}$ in
the argument above. Summing over the possible values of $u$ and
$v$ gives the two series.

We next compare $x$ with $\mu$. Uniformly in $\ell$, we have
$x/\mu=1+O(a+k/n)$. Integrating the derivative estimate
in Lemma~\ref{lem:series}, as in the preceding proof, gives
\[
 \left|\log\Phi_q(x)-\log\Phi_q(\mu)\right|
 =O_\eta\left(\frac{(a+k/n)\log L}{q}\right)
 =O_\eta\left(\frac{L^2\log L}{np^2}\right)=o(1).
\]
\indent By Lemma~\ref{lem:poissonization}, it follows that $\rho=(1+o(1))(1-a)^{n-k}\cdot \Phi_q(x)$, uniformly in $\ell$.

We proceed with the first inequality. Since
$\theta=a_{\ }e^{q\ell}\geq a$, we have
$b=(1-a)^2+a(\theta-a)$. Using $\log(1+w)\leq w$, $(1-a)^{-2}\leq 2$ and
$-\log(1-a)\leq2a$ for all large enough $n$, we obtain
 \[\log\left(\frac{b^{n-2k+\ell}}{(1-a)^{2(n-k)}}\right)
 =(n-2k+\ell)\log\left(1+\frac{a(\theta-a)}{(1-a)^2}\right)
 -2(k-\ell)\log(1-a)
 \leq2na(\theta-a)+4ka.\]
\indent Besides, $\Phi_q$ is increasing and $\Phi_q(w)\leq e^w$
for $w\geq0$, which helps us deduce that
\[
 J(x,\theta)\leq\Phi_q(x)^2\cdot \Phi_q(x\theta)
             \leq\Phi_q(x)^2\cdot e^{x\theta}.
\]
\indent Plugging in the estimates for $b$ and $J(x,\theta)$ into the joint-probability bound, we get:  
\[
 \bP(S,T\text{ both isolate})
 \leq(1-a)^{2(n-k)}\cdot \Phi_q(x)^2
      \exp\bigl(4ka+2na(\theta-a)+x\theta\bigr).
\]
\indent Recall that $\rho=(1+o(1))(1-a)^{n-k}\cdot \Phi_q(x)$.
Moreover, $x\sim\mu$, $na\sim\mu$ and $ka=o(1)$, hence
$4ka+2na(\theta-a)+x\theta\leq(3+o(1))\mu\theta+o(1)$. This proves the first inequality, as it implies that 
 $$\bP(S,T\text{ both isolate})
 \leq\rho^2\exp\bigl((3+o(1))\mu\theta+o(1)\bigr)
 \leq\rho^2\exp\bigl(4\mu\theta+o(1)\bigr).$$
\indent We now turn to the second inequality. We first compare the factor $b^{n-2k+\ell}$ with $(1-a)^{n-k}$, which appears in our expression for $\rho$. The identity $b=(1-a)-a(1-\theta)$ gives
\begin{align*}
 \log\left(\frac{b^{n-2k+\ell}}{(1-a)^{n-k}}\right)
 &=(n-2k+\ell)
   \log\left(1-\frac{a(1-\theta)}{1-a}\right)
   -(k-\ell)\log(1-a)\\
 &\leq-(n-2k+\ell)a(1-\theta)+2ka\leq-\mu(1-\theta)+3ka.
\end{align*}
\indent Here we wrote $n-2k+\ell=(n-k)-(k-\ell)$ and used $\log(1-w)\leq-w$,
$-\log(1-a)\leq2a$. Since $ka=o(1)$, the joint-probability bound therefore gives
\[
 \bP(S,T\text{ both isolate})
 \leq(1-a)^{n-k}J(x,\theta)\cdot
      \exp\bigl(-\mu(1-\theta)+o(1)\bigr).
\]
\indent Recall that $\rho=(1+o(1))(1-a)^{n-k}\cdot \Phi_q(x)$.
Thus we need to compare $J(x,\theta)$ with one copy of
$\Phi_q(x)$. Put $y:=x(1-\theta)$, so that the two
exclusive-group factors in $J(x,\theta)$ are
$\Phi_q(ye^{-qi})$. Since $S\ne T$, we have
$1-\theta\geq p$, which implies that $y\sim\mu(1-\theta)$ and $y\geq px=\Theta_\eta(L^2)\to\infty$.
It is therefore enough to prove that $J(x,\theta)\leq\Phi_q(x)\exp(o(y))$.

To begin, Lemma~\ref{lem:series} gives
$\log(qx)=2\log L+O_\eta(1)\to\infty$ and $\log\Phi_q(x)\sim \log^2(qx)/(2q)$.\\
\indent We now distinguish two cases. Suppose first that $qy\geq \log(qx)$. Then $qy\to\infty$, so
Lemma~\ref{lem:series} gives $\log\Phi_q(y)=O_\eta\left(q^{-1}\cdot\log^2(qy)\right)=o(y)$.
Using $\Phi_q(ye^{-qi})\leq\Phi_q(y)$ in the definition of $J(x,\theta)$, we obtain
$J(x,\theta) \leq\Phi_q(x\theta)\cdot \Phi_q(y)^2\leq\Phi_q(x)\exp(o(y))$.\\
\indent Suppose next that $qy<\log(qx)$. Then $2y<2\log(qx)\cdot q^{-1}=o\bigl(\log\Phi_q(x)\bigr)$.
 In particular, $e^{2y}\leq\Phi_q(x)^{1/6}$ for all
large enough $n$. We then split the sum defining
$J(x,\theta)$ as follows.

 For all $i\leq \log(qx)/(2q)$, we use $(x\theta)^i\leq x^i$ and
$\Phi_q(ye^{-qi})\leq e^y$. The lower-tail estimate in
Lemma~\ref{lem:series} bounds this contribution by
\[ e^{2y}\sum_{i\leq \log(qx)/(2q)}
       \frac{x^i}{i!}\cdot e^{-q\binom i2}
 \leq e^{2y}\cdot \Phi_q(x)^{5/6} \leq\Phi_q(x).
\]
\indent For $i>\log(qx)/(2q)$, we observe that $ye^{-qi}\leq ye^{-\log(qx)/2}=y/\sqrt{qx}$.
Using $\Phi_q(w)\leq e^w$, the contribution of these
terms is therefore at most
\[
 \exp\left(\frac{2y}{\sqrt{qx}}\right)
 \sum_{i> \log(qx)/(2q)}\frac{(x\theta)^i}{i!}e^{-q\binom i2}
 \leq\Phi_q(x)\exp\left(\frac{2y}{\sqrt{qx}}\right).
\]
\indent Adding the two contributions gives $J(x,\theta)
 \leq2\Phi_q(x)\exp\left(2y/\sqrt{qx}\right) =\Phi_q(x)\exp(o(y))$,
because $qx\to\infty$ and $y\to\infty$. This proves the required estimate in both cases.
Substituting it into the joint-probability bound, we obtain
\[
 \bP(S,T\text{ both isolate})
 \leq(1-a)^{n-k}\cdot \Phi_q(x)
      \exp\bigl(-\mu(1-\theta)+o(y)\bigr).
\]
\indent Plugging in $\rho=(1+o(1))(1-a)^{n-k}\cdot \Phi_q(x)$ and using
$y\sim\mu(1-\theta)\to\infty$, we conclude that
 $\bP(S,T\text{ both isolate}) \leq\rho\exp\bigl(-(1-o(1))\mu(1-\theta)\bigr)$.
Since $\theta=e^{-q(k-\ell)}$, this proves the second
inequality.

Finally, suppose that $q\ell\leq1$, so that $\theta\leq ea$.
For the sharper estimate, we fix the exclusive-group
sizes $u,v$ and sum over the common-group size.
Rearranging the same series gives
\[
 J(x,\theta)=\sum_{u,v\geq0}
 \frac{\bigl(x(1-\theta)\bigr)^{u+v}}{u!\cdot v!}\cdot
 e^{-q\binom u2-q\binom v2}\cdot 
 \Phi_q\bigl(x\theta e^{-q(u+v)}\bigr).
\]
\indent Recall that $u$ counts the vertices outside $S\cup T$
that are undominated by $S$ but dominated by $T$, while
$v$ counts those undominated by $T$ but dominated by $S$. If both $u$ and $v$ are greater than $\log(qx)/(2q)$, then $x\theta e^{-q(u+v)}\leq\theta/q$.
These terms therefore contribute at most $\Phi_q(x)^2e^{\theta/q}$.
For the other terms, we once again apply the lower-tail
estimate in Lemma~\ref{lem:series}. 
First consider the terms with $2qu\leq\log(qx)$. Since
$x(1-\theta)\leq x$, Lemma~\ref{lem:series} gives
\[
 \sum_{u=0}^{\lfloor\log(qx)/(2q)\rfloor}
 \frac{\bigl(x(1-\theta)\bigr)^u}{u!}\cdot e^{-q\binom u2}
 \leq
 \sum_{u=0}^{\lfloor\log(qx)/(2q)\rfloor}\frac{x^u}{u!}e^{-q\binom u2}
 \leq\Phi_q(x)^{5/6}.
\]
\indent There is no restriction on $v$, so its corresponding sum is
\[
 \sum_{v\geq0}
 \frac{\bigl(x(1-\theta)\bigr)^v}{v!}\cdot e^{-q\binom v2}
 =\Phi_q\bigl(x(1-\theta)\bigr)\leq\Phi_q(x).
\]
\indent Lastly, the factor accounting for vertices which are left undominated
by both sets satisfies
$\Phi_q(x\theta e^{-q(u+v)})\leq e^{x\theta}$.
Thus all terms with $2qu\leq\log(qx)$ contribute at most $e^{x\theta}\cdot \Phi_q(x)^{11/6}$.
By symmetry, the same bound holds for the terms with
$2qv\leq\log(qx)$. Their combined contribution is therefore
at most $2e^{x\theta}\cdot \Phi_q(x)^{11/6}$.

Our hypothesis on $p$ rewrites $np\gg\sqrt n\log^{3/2}n$, which implies $\theta/q=O_\eta\left(L^2/np^2\right)=o(1)$ and $na^2q=O_\eta\left(L^4/np\right)=o(1)$.
Since $x\sim na$, $\theta\leq ea$ and
$\log\Phi_q(x)\sim \log^2(qx)/(2q)$, we also have $x\theta\cdot \log^{-1}\Phi_q(x)=O_\eta\left(na^2q\cdot \log^{-2}(qx)^2\right)=o(1)$.
It follows that
\[
 \frac{J(x,\theta)}{\Phi_q(x)^2}
 \leq e^{\theta/q}
     +2\exp\left(x\theta-\frac16\log\Phi_q(x)\right)
 =1+o(1).
\]
\indent Moreover, $\theta-a=a(e^{q\ell}-1)=O(aq\ell)$. Since $ka=o(1)$ and $na^2q=o(1)$, we obtain
\[
 \log\left(\frac{b^{n-2k+\ell}}{(1-a)^{2(n-k)}}\right)
 \leq4ka+2na(\theta-a)
 =4ka+O(na^2q\ell)=o(\ell+1).
\]
\indent Combining this with the estimate for $J(x,\theta)$
and our expression for $\rho$, we conclude that
$\bP(S,T\text{ both isolate})
 \leq\rho^2\exp\bigl(o(\ell+1)\bigr)$,
completing the proof.
\end{proof}

Our next goal is to estimate the second moment
of $Y_k$. Most pairs of $k$-sets have a small intersection,
and the last bound in Lemma \ref{lem:pair} will control their contribution.
For pairs with a larger intersection, we compare their
number with the probability that both sets are isolating.

\begin{prop}\label{prop:second-moment}
Fix $\eta\in(0,1)$ and let $G\sim G(n,p)$, where
$0<p\leq1-\eta$ and $np^2/\log^3 n\to\infty$.
Let $k=k(n)$ be an integer satisfying $qk=\log(np)-2\log\log(np)+O(1)$. Then
\begin{equation*}\label{eq:second-moment}
 \bE[Y_k^2]\leq(1+o(1))M_k^2+(1+o(1))M_k.
\end{equation*}
In particular, if $M_k\to\infty$, then $G$ has an isolating
set of size $k$ whp.
\end{prop}

\begin{proof}
Recall that $L:=\log(np)$. As in Lemma \ref{lem:pair}, fix a $k$-set $S$, write $\mu:=(n-k)(1-p)^k$ and
let $\rho:=\bP(S\text{ isolates})$. By symmetry, $M_k=\binom nk\rho$.
Our assumptions imply $qk\sim L=\Theta(\log n)$, $q\mu=\Theta_\eta(L^2)$ and $k^2/n=O_\eta\bigl(L^2/(np^2)\bigr)=o(1)$.

The random variable $Y_k^2$ counts the ordered pairs $(S,T)$
of isolating $k$-sets, allowing $S=T$. Writing $\mathcal{S}_k$ for the family of vertex $k$-sets, we get $\bE[Y_k^2]
 =\sum_{S,T\in \mathcal{S}_k}\bP(S,T\text{ both isolate})$.
The pairs with $S=T$ contribute exactly $M_k$.
We estimate the other pairs according to
$\ell:=|S\cap T|$.
For each fixed $S$, there are
$\binom k\ell\binom{n-k}{k-\ell}$ choices of $T$
with this intersection size.

We first consider pairs with
$0<q(k-\ell)<3\log L$, which rewrites as $k-3q^{-1}\log L\leq l\leq k-1$.
The second bound in Lemma~\ref{lem:pair} shows that their
total contribution is at most
\[
 M_k
 \sum_{\ell=k-\lfloor3q^{-1}\log L\rfloor}^{k-1}
 \binom k\ell\binom{n-k}{k-\ell}
 \exp\bigl(-(1-o(1))\mu(1-e^{-q(k-\ell)})\bigr).
\]
\indent We show that the sum multiplying $M_k$ is of order $o(1)$.
If $q(k-\ell)\leq1$, then $1-e^{-q(k-\ell)}\geq q(k-\ell)/2$.
On the other hand, the binomial coefficient bounds give
\[
 \log\left(\binom k\ell\binom{n-k}{k-\ell}\right)
 \leq(k-\ell)
      \log\left(\frac{e^2kn}{(k-\ell)^2}\right)
 =O\bigl((k-\ell)L\bigr).
\]
\indent Since $q\mu=\Theta_\eta(L^2)$, the negative exponent
outweighs this number of choices. Thus, for some fixed
$c>0$, the sum of these terms is at most
$\sum_{\ell} \exp\bigl(-cL^2(k-\ell)\bigr)
 \leq e^{-cL^2}/\bigl(1-e^{-cL^2}\bigr)=o(1)$.
Conversely, when $1<q(k-\ell)<3\log L$, the exponential factor is
at most $e^{-\mu/3}$ for all large enough $n$.
Since there are at most $1+3\log L/q$ possible intersection
sizes and each product of binomial coefficients is
at most $(kn)^{3\log L/q}$, the sum of the terms with $1<q(k-\ell)<3\log L$ is at most
$\exp\left(-\dfrac{\mu}{3}+O_\eta\left(\dfrac{\log L\log n}{q}\right) \right)=o(1)$,
where we used $\dfrac{\log L\log n}{q\mu}=O_\eta\left(\dfrac{\log L}{L}\right)=o(1)$.
We conclude that all pairs with
$0<q(k-\ell)<3\log L$ have contribution $o(M_k)$.

We now consider the pairs for which
$q(k-\ell)\geq3\log L$, or equivalently $l<k-3q^{-1}\log L$.
When $\ell=0$, the last bound in Lemma~\ref{lem:pair}
gives joint probability at most $(1+o(1))\rho^2$.
Since there are at most $\binom nk^2$ ordered pairs,
their contribution is at most $(1+o(1))M_k^2$.
For $\ell\geq 1$, we proceed by estimating the proportion $P(k,l)$
of $k$-sets meeting a fixed $k$-set in exactly $\ell$
vertices. Note that a uniformly chosen $k$-set contains any prescribed
$\ell$-subset with probability
$\prod_{i=0}^{\ell-1}\frac{k-i}{n-i}\leq(k/n)^\ell$. Thus, a union bound over the $\ell$-subsets of $S$ gives
\[
 P(k,l)\leq\frac{\binom k\ell\binom{n-k}{k-\ell}}{\binom nk}
 \leq\binom k\ell\left(\frac{k}{n}\right)^\ell
 \leq\left(\frac{ek^2}{n\ell}\right)^\ell.
\]
\indent Since $M_k^2=\binom nk^2\rho^2$, we can multiply this
proportion by a bound for the joint probability divided
by $\rho^2$ to estimate the contribution divided by $M_k^2$.
Suppose first that $0<q\ell\leq1$.
By the last bound in Lemma~\ref{lem:pair}, the contribution
for each such intersection size, divided by $M_k^2$,
is at most $(ek^2/nl)^l \exp\bigl(o(\ell+1)\bigr)$.
For large enough $n$, the error in the exponent
is at $\leq \ell+1$.
Summing over these intersection sizes gives the upper bound
\[
 e\sum_{\ell\geq1}
       \left(\frac{e^2k^2}{n\ell}\right)^\ell
 \leq e\sum_{\ell\geq1}
       \left(\frac{e^2k^2}{n}\right)^\ell
 =O\left(\frac{k^2}{n}\right)=o(1).
\]
\indent Finally, suppose that $q\ell>1$ and
$q(k-\ell)\geq3\log L$.
The first bound in Lemma~\ref{lem:pair}, together with the estimate on $P(k,l)$, bounds the contribution for
each such $\ell$, divided by $M_k^2$, by
\[
 \exp\left(
 -\ell\log\left(\frac{n\ell}{ek^2}\right)
 +4\mu e^{-q(k-\ell)}+o(1)
 \right).
\]
\indent We now compare the two main terms in this exponent.
Since $\ell>1/q$ and $qk\sim L$, we have $\log\left(n\ell/ek^2\right)\geq\log n -\log(eqk^2)=L-2\log L+O_\eta(1)$. Thus the negative term in exponential above is at most $-(1-o(1))\ell L$.
For the positive term, observe that $e^{qu}/u$ is
increasing if $qu>1$.\\
The condition $q(k-\ell)\geq3\log L$ implies
$\ell\leq k-3\log L/q$. Hence, using $qk\sim L$ and $q\mu=\Theta_\eta(L^2)$, which also implies that $k-3\log L/q\sim k$, we get that 
\[
 \frac{\mu e^{-q(k-\ell)}}{\ell L}
 =\frac{\mu e^{-qk}}{L}\cdot\frac{e^{q\ell}}{\ell}
 \leq\frac{\mu e^{-qk}}{L}\cdot
       \frac{e^{q(k-3\log L/q)}}{k-3\log L/q}
 =\frac{\mu L^{-3}}{(k-3\log L/q)L}
 =O_\eta(L^{-3})=o(1).
\]
\indent Consequently, the contribution for each such $\ell$,
divided by $M_k^2$, is at most $e^{-\ell L/2}$ for
all large enough $n$. Summing up gives $\sum_{\ell\geq1}e^{-\ell L/2}=o(1)$. 
Adding the diagonal contribution and all the other
contributions, we obtain $\bE[Y_k^2]
 \leq M_k+o(M_k)+(1+o(1))M_k^2+o(M_k^2)$,
which proves the desired estimate.
If $M_k\to\infty$, this implies that
$\operatorname{Var}(Y_k)=o(M_k^2)$.
By Markov's inequality
we get $\bP(Y_k=0)
 \leq\bP\bigl((Y_k-M_k)^2\geq M_k^2\bigr)
 \leq\operatorname{Var}(Y_k)\cdot M_k^{-2}
 =o(1)$, completing the proof.
\end{proof}

We are now in a position to give a proof of
Theorem~\ref{thm:variable-counting}.
The expectation criterion follows from Markov's inequality
and Proposition~\ref{prop:second-moment}.
For the two-point concentration, we first locate the integer
$r$ appearing in the statement. We then compare the numbers
of isolating sets of consecutive sizes by counting their
extensions.

\begin{proof}[Proof of Theorem~\ref{thm:variable-counting}]
To begin with, note that every vertex set containing an isolating set is also isolating.
Consequently, $\iota(G)\leq k$ if and only if $Y_k>0$. On one hand, if $M_k\to0$, Markov's inequality gives $\bP(\iota(G)\leq k)=\bP(Y_k>0)\leq M_k=o(1)$.
On the other hand, when $M_k\to\infty$, Proposition~\ref{prop:second-moment}
gives $\bP(Y_k>0)=1-o(1)$. This proves the first assertion.

\indent Before turning to the last part,
we record a comparison between consecutive expectations.
Every isolating $k$-set has exactly $n-k$ extensions
of size $k+1$, all of which are isolating.
Conversely, an isolating $(k+1)$-set contains at most
$k+1$ isolating $k$-sets. By double counting these inclusions and
taking expectations, we obtain that $(n-k)Y_k\leq(k+1)Y_{k+1}$ and that $(n-k)M_k\leq(k+1)M_{k+1}$.
In particular, $M_k$ is non-decreasing for $k<n/2$.\\
\indent We next show that the smallest integer $r$ satisfying
$M_r\geq1/\log n$ lies in the range considered in
Proposition~\ref{prop:second-moment}.
Suppose first that
$qk=L-2\log L+O(1)$, and put
$\mu:=(n-k)(1-p)^k$.
Since $p\leq q\leq(1-\log\eta)p$, we have $\log\left(nk^{-1}\right) =L+\log\left(qp^{-1}\right)-\log(qk)
 =L-\log L+O_\eta(1)$. The binomial coefficient bounds therefore give
$\log\binom nk =k\log\left(nk^{-1}\right)+O(k)=(1+o(1))L^2/q$. Moreover, $(1-p)^k=o(1)$ and
$\mu=\Theta_\eta(L^2/q)$, so $(n-k)\log\bigl(1-(1-p)^k\bigr) =-\mu+o\left(L^2/q\right)$.
 Applying Lemma~\ref{lem:poissonization}, and using $\log\Phi_q(\mu)=o(L^2/q)$, we obtain that
\begin{align*}
 \log M_k
 &=\log\binom nk
   +(n-k)\log\bigl(1-(1-p)^k\bigr)
   +\log\Phi_q(\mu)+o(1)\\
 &=(1+o(1))\frac{L^2}{q}-(n-k)(1-p)^k.
\end{align*}
\indent Choose a sufficiently large constant $C>0$,
depending only on $\eta$, and set
\[
 k_-:=\left\lfloor\frac{L-2\log L-C}{q}\right\rfloor,
 \qquad
 k_+:=\left\lceil\frac{L-2\log L+C}{q}\right\rceil.
\]
\indent Since $q\leq-\log\eta$, both integers satisfy
$qk_\pm=L-2\log L+O_\eta(1)$. We also have
\[
 (n-k_-)(1-p)^{k_-}
 \geq(1-o(1))\cdot e^C L^2p^{-1},\qquad \text{and}
 \qquad
 (n-k_+)(1-p)^{k_+}
 \leq e^{-C}L^2p^{-1}.
\]
\indent As $1\leq q/p\leq1-\log\eta$, choosing $C$
large enough makes the first quantity larger than
$2L^2/q$ and the second less than $L^2/(2q)$.
Our estimate for $\log M_k$ then gives, for all large
enough $n$, that $ M_{k_-}\leq\exp\left(-L^2/2q\right)=o\left(\log^{-1} n\right)$ and $M_{k_+}\geq\exp\left(L^2/4q\right)\longrightarrow\infty$.
Since $k_+=o(n)$, the monotonicity of $M_k$ implies that $k_-<r\leq k_+$.
Consequently, $qr=L-2\log L+O_\eta(1)$.

\indent By the definition of $r$, we have $M_{r-1}<1/\log n$. Markov's inequality therefore shows that
\[
 \bP\bigl(\iota(G)<r\bigr)
 =\bP(Y_{r-1}>0)
 \leq M_{r-1}=o(1).
\]
\indent On the other hand, the estimate from above gives $M_{r+1}\geq\dfrac{n-r}{r+1}M_r
          \geq\dfrac{n-r}{(r+1)\log n}
          \longrightarrow\infty$.\\
\indent Indeed, $r=O_\eta(L/p)=o(n)$, while our
assumption on $p$ implies $np/\log^2n\to\infty$.
Finally, $q$ is bounded, so $r+1$ also satisfies
$q(r+1)=L-2\log L+O_\eta(1)$.
Proposition~\ref{prop:second-moment} applied at $r+1$
shows that $\iota(G)\leq r+1$ whp.
Together with the lower bound, this yields
$\iota(G)\in\{r,r+1\}$ whp, completing the proof.
\end{proof}

\section{The gap between isolation and domination}\label{sec:comparison}

In this section we prove Theorem~\ref{thm:gap} by comparing the expected
numbers of dominating and isolating sets of the same size.
Lemma~\ref{lem:poissonization} shows that an independent
set of undominated vertices increases the expectation by a factor
given by the isolation potential. We will compare this gain with
the change in the expectation when the size of the set increases.

Throughout this part, $\eta\in(0,1)$ is fixed and we assume, as in the hypothesis of
Theorem~\ref{thm:gap}, that $0<p=p(n)\leq1-\eta$ with $-\log p=o(\log n)$. We still set $q:=-\log(1-p)$ and
$L:=\log(np)$. In particular, $L\sim\log n$, and the assumptions
of Section~\ref{sec:concentration} hold. Let $X_k$ count the
dominating sets of size $k$ in $G(n,p)$; as before, $Y_k$ counts
the isolating sets and $M_k=\bE[Y_k]$.

We start with an existence criterion for dominating sets.
This follows from the same counting argument as
Proposition~\ref{prop:second-moment}, once we check the
corresponding joint probabilities.

\begin{lem}\label{lem:domination-existence}
Fix $\eta\in(0,1)$ and let $G\sim G(n,p)$, where
$0<p\leq1-\eta$ and $np^2/\log^3 n\to\infty$.
Let $X_k$ count the dominating sets of size $k=k(n)$ in $G$,
where $kq=\log(np)-2\log\log(np)+O(1)$.
If $\bE[X_k]\to\infty$, then whp $G$ has a dominating set of size
$k$.
\end{lem}

\begin{proof}
Set $a:=(1-p)^k$ and let $S$ and $T$ be two distinct $k$-sets. Write $\ell:=|S\cap T|$.
For a vertex $v$ outside $S\cup T$, the probability that
$v$ has no neighbor in $S$ is $a$, and the same holds for
$T$. Moreover, $v$ has no neighbor in $S\cup T$ precisely
when all its $2k-\ell$ possible edges to this union are
absent, which happens with probability $(1-p)^{2k-\ell}$.
By inclusion-exclusion, the probability that $v$ has
a neighbor in each of $S$ and $T$ is therefore $b:=1-a-a+(1-p)^{2k-\ell}
   =1-2a+a(1-p)^{k-\ell}$.
   
\indent These events are independent for different vertices
outside $S\cup T$. Ignoring the domination requirements
inside $S\cup T$ therefore gives $\bP(S,T\text{ both dominate})\leq b^{n-2k+\ell}$.
Recall also that $\bP(S\text{ dominates})=(1-a)^{n-k}$.
The comparisons used in the proof of
Lemma~\ref{lem:pair} give
\begin{align*}
 \log\left(\frac{b^{n-2k+\ell}}{(1-a)^{2(n-k)}}\right)
 &\leq4ka+2na\bigl((1-p)^{k-\ell}-a\bigr),\\
 \log\left(\frac{b^{n-2k+\ell}}{(1-a)^{n-k}}\right)
 &\leq3ka-(n-k)a\bigl(1-(1-p)^{k-\ell}\bigr).
\end{align*}
\indent As in that proof, our assumptions imply
$k=o(n)$, $ka=o(1)$, $p(n-k)a=\Theta_\eta(\log^2(np))$ and $na^2q=o(1)$.
Since $S\ne T$, we have
$1-(1-p)^{k-\ell}\geq p$.
Consequently, these comparisons imply the first two
bounds in Lemma~\ref{lem:pair}, with domination
in place of isolation.

For the third bound, suppose that $q\ell\leq1$.
Then
$(1-p)^{k-\ell}-a=a(e^{q\ell}-1)=O(aq\ell)$,
so the first comparison gives $\log\bigl(b^{n-2k+\ell}\cdot (1-a)^{-2(n-k)}\bigr)
 \leq4ka+O(na^2q\ell)=o(\ell+1)$.
Thus all three joint-probability bounds used in
Proposition~\ref{prop:second-moment} hold for domination.
Counting pairs according to their intersection size
exactly as in that proposition, we obtain
\[
 \bE[X_k^2]
 \leq(1+o(1))\bE[X_k]^2+(1+o(1))\bE[X_k].
\]
\indent Since $\bE[X_k]\to\infty$, it follows that
$\operatorname{Var}(X_k)=o(\bE[X_k]^2)$, thus $\bP(X_k=0)
 \leq\dfrac{\operatorname{Var}(X_k)}{\bE[X_k]^2}
 =o(1)$, 
completing the proof.
\end{proof}

The next step is to compare continuous approximations to the logarithms
of the two expectations. Their zeros locate the graph
parameters up to a bounded error. We estimate the distance
between these zeros using the isolation potential and the
derivatives of the two functions.

\begin{lem}\label{lem:centres}
Fix $\eta\in(0,1)$ and let $p=p(n)$ satisfy
$0<p\leq1-\eta$ and $\log(1/p)=o(\log n)$.
Define $f(u):=u\log\left(enu^{-1}\right)-(n-u)e^{-qu}$ on the interval $(0,n)$. 
For all sufficiently large $n$, the functions $f(u)$ and
$f(u)+\log\Phi_q\bigl((n-u)e^{-qu}\bigr)$ have unique zeros
$\kappa_\gamma$ and $\kappa_\iota$, respectively, both in the closed interval with endpoints  $q^{-1}(\log(nq)-2\log L-1)$ and  $q^{-1}(\log(nq)-2\log L+1)$.
These zeros satisfy $\kappa_\gamma-\kappa_\iota\sim 2q^{-1}(L^{-1}\log L)^2$.
Moreover, whp
\[
 \gamma(G(n,p))=\kappa_\gamma+O(1)
 \qquad\text{and}\qquad
 \iota(G(n,p))=\kappa_\iota+O(1).
\]
\end{lem}

\begin{proof}
Write $\mu(u):=(n-u)e^{-qu}$ for convenience. We first show how these two functions relate to the expected numbers of dominating and isolating sets. If
$qk=L-2\log L+O(1)$ for integers $k$, then $k=n^{o(1)}$ and
$e^{-qk}=n^{-1+o(1)}$. Comparing $\log(k!)$ with the integral
of $\log x$ gives $\log(k!)=k\log k-k+O(\log k)$, and hence
\[
 \log\binom nk
 =k\log\left(\frac{en}{k}\right)
   +O\left(\log k+\frac{k^2}{n}\right)
 =k\log\left(\frac{en}{k}\right)+O(\log k).
\]
\indent Also,
$(n-k)\log(1-e^{-qk})=-\mu(k)+O(ne^{-2qk})=-\mu(k)+o(1)$.
 Lemma~\ref{lem:poissonization} then gives
\begin{equation}\label{eq:curves-accurate}
\begin{aligned}
 \log\bE[X_k]&=f(k)+O(\log k),\\
 \log M_k&=f(k)+\log\Phi_q(\mu(k))+O(\log k).
\end{aligned}
\end{equation}
\indent We now locate the zeros. Since $q=\Theta_\eta(p)$, any point $u$ from the
interval in the statement satisfies $qu=L-2\log L+O(1)$. In this range,
$\log(n/u)=\log(nq)-\log(qu)\sim L$ since $u\sim L/q=o(n)$, so $u\log(en/u)\sim L^2/q$. At the left and right endpoints,
$q\mu(u)$ is asymptotic to $eL^2$ and $e^{-1}L^2$,
respectively. Thus $f$ is negative at the left endpoint and positive at the right endpoint. Moreover,
$q\mu(u)=\Theta(L^2)$ throughout the interval, so Lemma~\ref{lem:series} gives us that uniformly 
$\log\Phi_q(\mu(u)) \sim2q^{-1}\log^2 L =o\left(L^2/q\right)$. This implies that the second function has the same signs at the endpoints. Moreover, both these functions are strictly increasing. Indeed, $\mu'(u)<0$,
and differentiating the defining series gives
$\Phi_q'(x)=\Phi_q(xe^{-q})\leq\Phi_q(x)$.
It follows that $f'(u)=\log(n/u)-\mu'(u)\geq\log(n/u)>0$ and
\[
 \frac{d}{du}\bigl(f(u)+\log\Phi_q(\mu(u))\bigr)
 =\log\left(\frac{n}{u}\right)
  -\mu'(u)\left(1-\frac{\Phi_q'(\mu(u))}{\Phi_q(\mu(u))}\right)
 \geq\log\left(\frac{n}{u}\right)>0.
\]
\indent Thus each function has a unique zero in the interval given in the hypothesis. Since $\Phi_q(\mu(u))>1$, we have also get $\kappa_\iota<\kappa_\gamma$.
At either zero, the equation defining it gives
$q\mu(u)\sim L^2$: the term $u\log(en/u)$ is asymptotic
to $L^2/q$, and the contribution of $\log\Phi_q(\mu(u))$
at the isolation zero is smaller. Since $\mu$ is decreasing,
we also have $q\mu(u)\sim L^2$ uniformly between the zeros.
Using $nq\to\infty$ and $u=o(n)$, we obtain
$f'(u)=\log\left(n/u\right) +\bigl(q+(n-u)^{-1}\bigr)\mu(u) \sim L^2$ for all $u\in (\kappa_\iota,\kappa_\gamma)$. The mean value theorem therefore gives
\[
 (\kappa_\gamma-\kappa_\iota)L^2 \sim f(\kappa_\gamma)-f(\kappa_\iota)
 =\log\Phi_q(\mu(\kappa_\iota)) \sim2q^{-1}(\log L)^2.
\]

Finally, for either zero $\kappa$, consider the integers
$\lfloor\kappa\rfloor-1$ and $\lceil\kappa\rceil+1$.
Each is at least one unit from $\kappa$, and since $q$ is
bounded, the intervening points still satisfy
$qu=L-2\log L+O(1)$. Both derivatives are at least
$\log(n/u)\sim L$ in this range, whereas the errors
in~\eqref{eq:curves-accurate} are $O(\log k)=o(L)$.
Thus the expectation tends to zero at the
lower integer and to infinity at the upper one.

By Markov's inequality, whp $G(n,p)$ has no dominating set
of size $\lfloor\kappa_\gamma\rfloor-1$ and no isolating set
of size $\lfloor\kappa_\iota\rfloor-1$. This also rules out
smaller sets of either kind, since adding vertices to a set
preserves both domination and isolation.
On the other hand, Lemma~\ref{lem:domination-existence}
gives a dominating set of size $\lceil\kappa_\gamma\rceil+1$
whp, and Proposition~\ref{prop:second-moment} gives an
isolating set of size $\lceil\kappa_\iota\rceil+1$ whp.
Combining these bounds, we obtain the claimed $O(1)$ estimates.
\end{proof}

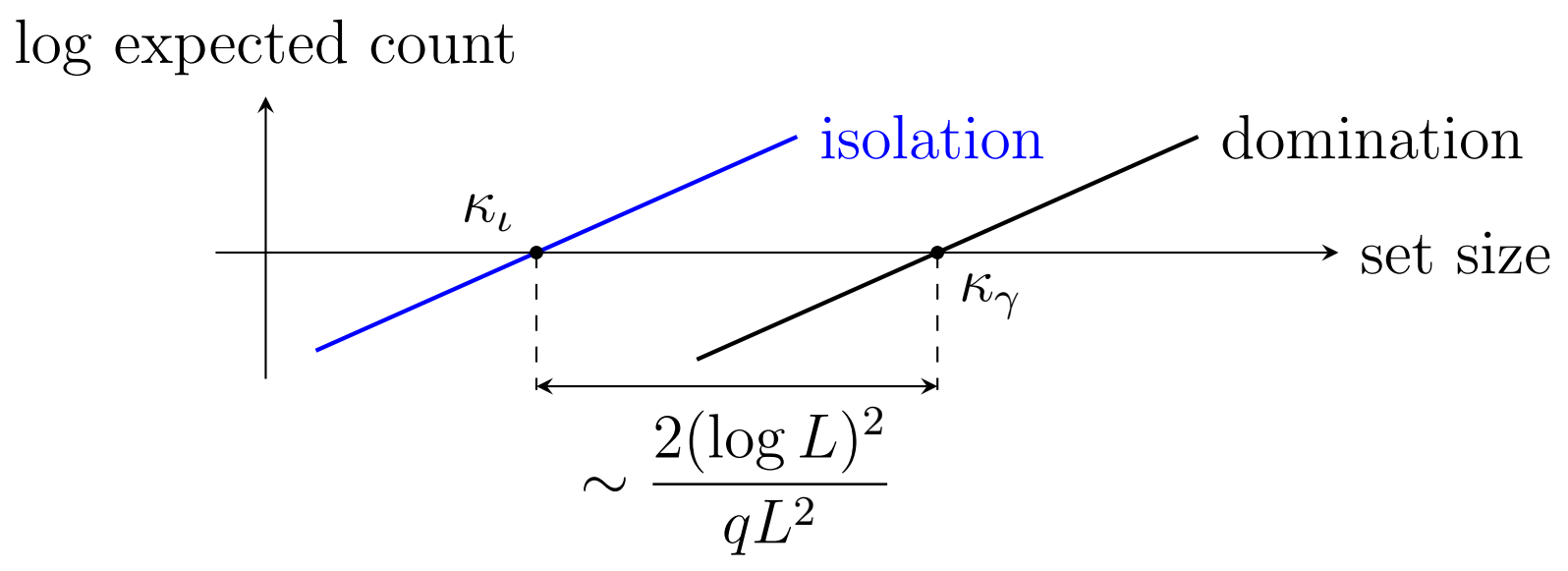
\begin{figure}[ht]
 \centering
 \begin{tikzpicture}[x=1.35cm,y=0.50cm,>=stealth]
  \draw[->] (0.4,0)--(6,0) node[right] {set size};
  \draw[->] (0.65,-1.7)--(0.65,2.1)
        node[above] {$\log$ expected count};
  \draw[thick] (2.8,-1.44)--(5.3,1.56)
        node[right] {domination};
  \draw[thick,blue] (0.9,-1.32)--(3.3,1.56)
        node[right] {isolation};
  \fill (2,0) circle (1.3pt);
  \fill (4,0) circle (1.3pt);
  \node[above left] at (2,0) {$\kappa_\iota$};
  \node[below right] at (4,0) {$\kappa_\gamma$};
  \draw[dashed] (2,0)--(2,-1.85);
  \draw[dashed] (4,0)--(4,-1.85);
  \draw[<->] (2,-1.8)--(4,-1.8);
  \node[below] at (3,-1.8)
        {$\displaystyle\sim\frac{2(\log L)^2}{qL^2}$};
 \end{tikzpicture}
 \caption{Schematic comparison of the two zeros and the
 common slope $\sim L^2$.}
 \label{fig:thresholds}
\end{figure}

\begin{proof}[Proof of Theorem~\ref{thm:gap}]
Recall $p_*:=(\log\log n/\log n)^2$. Since $L:=\log(pn)\sim\log n$, Lemma~\ref{lem:centres} gives
$\kappa_\gamma-\kappa_\iota\sim 2p_*/q$. The same lemma gives 
$\gamma(G(n,p))-\iota(G(n,p)) =\kappa_\gamma-\kappa_\iota+O(1)$ whp.\\
\indent If $p\ll p_*$, then $p=o(1)$, so $q\sim p$ and
$\kappa_\gamma-\kappa_\iota\sim 2p_*/p\to\infty$.
The bounded error is therefore negligible, proving~(i). If $p\geq cp_*$ for a fixed $c>0$, then $q\geq p$ gives $\kappa_\gamma-\kappa_\iota=O_c(1)$, which proves~(iii). Hence, it remains to prove~(ii). Suppose that $p\gg p_*$. Then $\kappa_\gamma-\kappa_\iota=o(1)$. Let
$r:=\lfloor\kappa_\gamma+1/2\rfloor$, so that $r$ is a
nearest integer to $\kappa_\gamma$. For all sufficiently
large $n$, we have
\[
 r-1\leq\kappa_\gamma-\frac12\leq\kappa_\iota-\frac14,
 \qquad\text{and}\qquad
 r+1\geq\kappa_\gamma+\frac12.
\]
\indent All points between these sizes and the two zeros
still satisfy $qu=L-2\log L+O(1)$, since $q$ is bounded. As shown in Lemma~\ref{lem:centres}, both
functions have derivative at least $\log(n/u)\sim L$ in this range. Thus, for large enough $n$, the derivative is at least $L/2$. Using the distances from the zeros and the $o(L)$ errors in~\eqref{eq:curves-accurate}, we obtain
$\log M_{r-1}\leq-L/8+o(L)\longrightarrow-\infty$ and $\log\bE[X_{r+1}]\geq\frac{L}{4}+o(L)\longrightarrow\infty$.
By Markov's inequality, whp there is no isolating
set of size $r-1$, and hence none of smaller size.
On the other hand, Lemma~\ref{lem:domination-existence}
gives a dominating set of size $r+1$ whp. Therefore $r\leq\iota(G(n,p))\leq\gamma(G(n,p))\leq r+1$ whp,
which proves~(ii).
\end{proof}

For probabilities of the form $p=(\log n)^{-C}$,
Theorem~\ref{thm:gap} takes the following form.

\begin{cor}\label{cor:polylog}
Let $G\sim G(n,p)$, where $p=(\log n)^{-C}$ for a fixed
constant $C>0$. If $C\geq2$, then whp
 $\gamma(G)-\iota(G)\sim 2(\log n)^{C-2}(\log\log n)^2$.
If $0<C<2$, then whp $0\leq\gamma(G)-\iota(G)\leq1$.
\end{cor}

\begin{proof}
Since $p/p_* =(\log n)^{2-C} (\log\log n)^{-2}$,
we observe that $p\ll p_*$ for $C\geq2$, while $p\gg p_*$
for $0<C<2$. In the first scenario, part~(i) of
Theorem~\ref{thm:gap} gives the claimed asymptotic since $2p_*/p=2(\log n)^{C-2}(\log\log n)^2$.
In the second one, the assertion follows from part~(ii).
\end{proof}

\section{Graphs of large minimum degree}\label{sec:regular}

We now prove Theorem~\ref{thm:extremal}. As in the proof of
Theorem~\ref{thm:sparse}, the idea is to show that every small
set leaves enough undominated vertices to contain a copy of $F$.
We will use the following estimate for an arbitrary fixed graph
$F$.

\begin{lem}\label{lem:large-sets-contain-f}
Let $F$ be a fixed nontrivial graph and let $0<p=p(n)<1$. If $t=t(n)\leq n$ is an integer satisfying
$tp^{e(F)}/\log n\to\infty$, then whp every $t$-vertex set in $G(n,p)$ spans a copy of $F$.
\end{lem}

\begin{proof}
If $F$ has no edges, the result is immediate for $t\geq|V(F)|$.
Otherwise, put $r:=|V(F)|$ and fix a set $T$ of $t$ vertices.
Choose a maximal family of $r$-subsets of $T$ whose pairwise
intersections have size at most one. By maximality, every
$r$-subset of $T$ shares a pair with some member of this family.
Each member shares a pair with at most
$\binom r2\binom{t-2}{r-2}$ of the $r$-subsets of $T$.
The number of elements in this maximal family is therefore at least
\[
 \frac{\binom tr}{\binom r2\binom{t-2}{r-2}}
 =\frac{\binom t2}{\binom r2^2}
\]
\indent On each member, fix a copy of $F$ and consider
the event that all its edges are present. These events use
disjoint sets of edges, so they are independent, and each has
probability $p^{e(F)}$. If $T$ contains no copy of $F$, then
none of these events occurs. Taking a union bound over all $T$ and using $tp^{e(F)}/\log n\to\infty$, the probability that some $t$-set contains no
copy of $F$ is at most
\[
 n^t\exp\left(-p^{e(F)}\frac{\binom t2}{\binom r2^2}\right)
 =o(1).
\]
\end{proof}

For the regular construction, we will use the following coupling theorem.

\begin{thm}[Gao, Isaev and McKay~{\cite{GIM}}]
\label{thm:regular-coupling}
Let $d=d(n)$ be an integer such that $nd$ is even and
$\min\{d,n-d\}\gg\log^4 n$. Then there exist probabilities
$p_-,p_+\sim d/n$, together with a coupling of $G_-\sim G(n,p_-)$,
$G_+\sim G(n,p_+)$ and a uniformly random $d$-regular graph $R$,
all on the same vertex set, such that whp $G_-\subseteq R\subseteq G_+$.
\end{thm}

\begin{proof}[Proof of Theorem~\ref{thm:extremal}]
We first deal with the upper bound. Given a graph $G$ with
$\delta(G)\geq d$, choose each vertex independently with probability
$p=\log(d+1)/(d+1)$, and then add all vertices left undominated.
A vertex $v$ is left undominated by the initial choice precisely when
neither $v$ nor any of its neighbors is chosen. Since the choices are
independent and $\deg_G(v)\geq d$, the probability of this event is
$(1-p)^{\deg_G(v)+1} \leq (1-p)^{d+1} \leq e^{-p(d+1)}=(d+1)^{-1}$.
Thus the expected size of the
resulting dominating set is at most
\[
 \left(p+\frac{1}{d+1}\right)|V(G)|
 =\frac{\log(d+1)+1}{d+1}|V(G)|.
\]
\indent As a result, there is a dominating set of at most this size. Since every dominating set is also $F$-isolating, this proves the upper bound for both assertions.

For the lower bounds, take $n$ to be the least even integer satisfying $n\geq d(\log d)^2$. We first consider any probability $p\sim d/n$. Then $p\sim(\log d)^{-2}$ and
$\log(np)\sim\log n\sim\log d$. Put
$t:=\lceil(\log d)^{2e(F)+2}\rceil$. Since
$tp^{e(F)}\sim(\log d)^2$, Lemma~\ref{lem:large-sets-contain-f}
shows that whp every $t$-set in $G(n,p)$ contains a copy of $F$.
We next estimate how many vertices a small set leaves
undominated. As $(np)^\eps=(\log(np))^{2e(F)+3}\geq e^4\log^2(np)$ for large enough degree $d$, by Lemma~\ref{lem:expansion} applied with
$\eps=(2e(F)+3)\log\log(np)/\log(np)$, we get whp that every
set of size at most
\[
 k(p):=\frac{\log(np)-(2e(F)+3)\log\log(np)}{-\log(1-p)}
\]
leaves at least $(\log(np))^{2e(F)+3}/(3p)>t$ vertices
undominated. By putting together these two estimates, we deduce that
$\iota(G(n,p),F)>k(p)$ whp.

To make the minimum degree at least $d$, take
$p=(d+d^{2/3})/n$. The expected degree of each vertex is
$(n-1)p=d+(1+o(1))d^{2/3}$. By Chernoff's inequality and a
union bound, the probability that some vertex has degree
less than $d$ is at most
\[
 n\exp\left(-\frac{((n-1)p-d)^2}{2(n-1)p}\right)
 =n\exp\left(-\frac{(d^{2/3}-p)^2}{2(d+d^{2/3}-p)}\right)
 =n\exp\left(-\bigl(\tfrac12+o(1)\bigr)\sqrt[3]{d}\right)
 =o(1).
\]
\indent For this choice of $p$, we have $\log(np)=\log d+o(1)$ and
$n(-\log(1-p))=d(1+O((\log d)^{-2}))$. Consequently, whp
$\delta(G(n,p))\geq d$ and
\[
 \frac{\iota(G(n,p),F)}{n}
 \geq \frac{\log d-(2e(F)+3)\log\log d-O_F(1)}{d}.
\]
\indent This proves the first assertion. If $F$ is a forest,
the lower bound in the proof of Theorem~\ref{thm:sparse}
allows us to replace $k(p)$ by
$(\log(np)-2\log\log(np)-4)/(-\log(1-p))$. With the same
choice of $p$, this gives the sharper lower bound
$(\log d-2\log\log d-O_F(1))/d$ for the ratio.

It remains to prove the lower bound for $d$-regular graphs. Since $nd$ is even and
$\min\{d,n-d\}=d\gg\log^4 n$, we may take $G_-,R,G_+$ as in
Theorem~\ref{thm:regular-coupling}.
Since $p_-,p_+\sim d/n$, the estimates above apply to the
coupled binomial graphs. Thus, whp the
inclusions $G_-\subseteq R\subseteq G_+$ hold, every $t$-set
in $G_-$ contains a copy of $F$, and every set of size at most
$k(p_+)$ leaves more than $t$ vertices undominated in $G_+$.
For any such set $S$, let $U_S$ be the vertices it leaves
undominated in $G_+$. Since $R\subseteq G_+$, these vertices
are also undominated by $S$ in $R$. As $|U_S|>t$, the graph
$G_-[U_S]$ contains a copy of $F$, which also lies in $R[U_S]$
because $G_-\subseteq R$. Thus $S$ is not $F$-isolating in $R$.\\
\indent We conclude that whp $\iota(R,F)>k(p_+)$. Moreover, since $np_+\sim d$
and $-\log(1-p_+)\sim p_+$, we have $k(p_+)/n\sim(\log d)/d$.
Together with the upper bound, this proves the second assertion.
\end{proof}

\section{Concluding remarks}\label{sec:concluding-remarks}

The first-order estimate of Theorem~\ref{thm:sparse} extends
to every fixed nontrivial graph $F$ when
$p=o(1)$ and $\log(1/p)=o(\log n)$. Indeed, for any fixed
$\eps\in(0,1)$, Lemma~\ref{lem:expansion} shows that every set
of at most $(1-\eps)\log(np)/(-\log(1-p))$ vertices leaves
at least $n^{\eps/2}$ vertices undominated whp.
Lemma~\ref{lem:large-sets-contain-f} guarantees that every
set of $\lceil n^{\eps/2}\rceil$ vertices contains a copy of
$F$ whp. Together with $\iota(G,F)\leq\gamma(G)$, this gives
the same first-order asymptotics as for domination.
The triangle example at the end of
Section~\ref{sec:first-order} shows that the behavior can
change when $p$ is smaller. For a fixed graph $F$ containing
a cycle, what is the asymptotic behavior of
$\iota(G(n,p),F)$ when $p=n^{-\alpha}$ with $0<\alpha<1$?

For the domination threshold, Bohman, Warnke and Zhu~\cite{BWZ} proved
two-point concentration down to $p=n^{-2/3}\log^3 n$. It is natural to ask whether
the isolation number also concentrate on two consecutive
values for $p=n^{-\alpha}$ with $1/2\leq\alpha<2/3$.
One difficulty is already visible at $p=n^{-1/2}$. Two
independent uniformly random $k$-sets with $k\sim\log(np)/p$
have expected intersection size $k^2/n$ of order $\log^2 n$.
The proof of Proposition~\ref{prop:second-moment} uses
$k^2/n=o(1)$ to show that the contribution from pairs with
small nonzero intersection is negligible. In the sparser
range, these pairs require a different estimate.  The bipartite estimates of
Bohman, Warnke and Zhu~\cite{BWZ} suggest a possible starting point,
but the edge restrictions inside the undominated sets must also be
controlled.

Regarding the gap between isolation and domination, at the transition $p\sim cp_*$ for a fixed $c>0$,
Lemma~\ref{lem:centres} gives $\kappa_\gamma-\kappa_\iota\to2/c$.
If both zeros stay a fixed positive distance from the integers,
the argument of that lemma applies at their floors and ceilings:
the expected count tends to zero at the floor and to infinity
at the ceiling of each zero. Consequently, whp
\[
 \gamma(G(n,p))-\iota(G(n,p))
 =\lceil\kappa_\gamma\rceil-\lceil\kappa_\iota\rceil.
\]
\indent The case in which a zero approaches an integer is more delicate.
Suppose that $k$ satisfies the hypotheses of
Theorem~\ref{thm:variable-counting} and that $M_k\to\lambda$
for some $\lambda\in(0,\infty)$. The comparison of consecutive
expectations in the proof gives $M_{k-1}\to0$ and
$M_{k+1}\to\infty$, so $\iota(G(n,p))\in\{k,k+1\}$ whp.
Can we determine the probabilities of these two values? In
particular, we suspect that $Y_k$ might converge in distribution to a Poisson
random variable with mean $\lambda$. This would then give
$\bP(\iota(G(n,p))=k)\to1-e^{-\lambda}$.

\section*{Origin story of the paper and the use of AI}

After I finished my PhD, I kept reading and learning more mathematics, broadening my expertise, but I still found it hard to decide what problems to work on. One thing I missed was talking through the possible choices with my supervisor, Eoin Long. He has a tremendous sense of which questions I would enjoy and where I might make progress. I appreciated that even more once I had to make these choices on my own.\\
\indent At some point, confused as I was, I figured I could try to ask ChatGPT. I prompted it about my research and, after a few rounds of suggestions, the isolation number came up. It reminded me of something Eoin had taught me: if you want to get a feel for a graph parameter, try studying it in a random graph model. There was already a lot of work on domination in random graphs, so I had something to compare isolation with. My guess was that the two numbers wouldn’t be too far apart. That was the starting point for this project.\\
\indent I kept using ChatGPT and Codex as I worked on the paper, mostly to write out formulas, check calculations, and revise the writing. They also helped me work out some technical parts of the proofs.
 Once I had a draft, I was curious to see what the model would do with it. Using D\"{o}m\"{o}t\"{o}r P\'{a}lv\"{o}lgyi’s access to ChatGPT Pro, I gave it the draft and asked it to build on the results while trying to emulate my style of writing. I’ve included its response in the appendix as it came back. I take responsibility for the paper as a whole, including the AI mathematics in the appendix.

\appendix
\section{Appendix with AI results}\label{app:further}

In this appendix there are several consequences of the arguments above, fully derived by ChatGPT-6 Astra. First, a density parameter gives a first-order estimate for every fixed forbidden graph, and also improves the lower bound in the minimum-degree problem. We then consider two contrasting regimes for cyclic forbidden
graphs. Finally, we describe the vertices left undominated by an
isolating set, refine the comparison with domination, and prove the
Poisson limit suggested in Section~\ref{sec:concluding-remarks}.

\subsection{A first-order estimate for every fixed forbidden graph}\label{app:general-F}

The forest assumption in Theorem~\ref{thm:sparse} can be replaced by a
parameter which measures the density of the forbidden graph. For a fixed
nonempty graph $F$, write $v(H)=|V(H)|$ and $e(H)=|E(H)|$, and put
\[
 \rho(F):=\max_{\substack{H\subseteq F\\v(H)\geq1}}
                 \frac{e(H)+1}{v(H)},
 \qquad [x]_+:=\max\{x,0\}.
\]
Here $H$ is a subgraph, not necessarily an induced subgraph. In particular,
$H=K_1$ is allowed, so $\rho(F)\geq1$. If $F$ is a forest, every nonempty
subgraph $H$ satisfies $e(H)\leq v(H)-1$, and hence $\rho(F)=1$.
Conversely, a cycle in $F$ gives $\rho(F)>1$.

\begin{thm}\label{app:thm-general-F}
Let $F$ be a fixed nonempty graph, and let $0<p=p(n)<1$ satisfy
$np\to\infty$. Write $q=-\log(1-p)$ and $L=\log(np)$. Then whp
\[
 \iota(G(n,p),F)
 =\frac{[\log(n p^{\rho(F)})]_+}{q}
   +O_F\left(1+\frac{\log L}{q}\right).
\]
More precisely, for a constant $C_F>0$, whp
\begin{equation}\label{app:eq-general-F-bounds}
 \frac{\log(n p^{\rho(F)})-2\log L-C_F}{q}
 <\iota(G(n,p),F)
 \leq
 \left\lceil\frac{[\log(n p^{\rho(F)})+\log L]_+}{q}\right\rceil
 +\frac1p.
\end{equation}
\end{thm}

Thus the additional $1$ in the numerator defining $\rho(F)$ has a
concrete role. After choosing an initial set, we can remove the remaining
copies of a suitable subgraph of $F$ at a cost of order $1/p$. The density
parameter records when the expected number of these copies reaches that
scale.
We use the following standard consequence of Janson's inequality; see,
for example, Theorem~5.4 in~\cite{FKaronski}, originating in~\cite{JLR1990}. We include the overlap
calculation to make its application explicit.

\begin{lem}\label{app:lem-janson-F}
For every fixed graph $F$ with at least one edge, there is a constant
$c_F>0$ such that, for every sufficiently large integer $t$ and every
$0<p<1$,
\[
 \bP\bigl(G(t,p)\text{ contains no copy of }F\bigr)
 \leq \exp\left(-c_F
       \min_{\substack{H\subseteq F\\e(H)\geq1}}
              t^{v(H)}p^{e(H)}\right).
\]
\end{lem}

\begin{proof}
Let $Z$ count the copies of $F$, and write $f=v(F)$, $m=e(F)$ and
$\mu=\bE [Z]$. Then $\mu=\Theta_F(t^f p^m)$. For two copies whose
common vertices and common edges form a graph $H$ with at least one edge,
their contribution to the dependency sum in Janson's inequality is at
most $O_F(t^{2f-v(H)}p^{2m-e(H)})$. Copies with no common edge use
disjoint edge variables and do not contribute. Consequently, writing
$\Delta$ for this dependency sum, we have
\[
 \mu+\Delta
 \leq C_F\mu^2
       \sum_{\substack{H\subseteq F\\e(H)\geq1}}
            \frac{1}{t^{v(H)}p^{e(H)}}.
\]
The term $H=F$ also accounts for $\mu$. There are only finitely many
terms, with the number depending on $F$. Janson's inequality in the form
$\bP(Z=0)\leq\exp(-\mu^2/(2(\mu+\Delta)))$ proves the claim after
adjusting the constant.
\end{proof}

\begin{proof}[Proof of Theorem~\ref{app:thm-general-F}]
Put $\rho=\rho(F)$ and $B=\log(n p^\rho)$.
For the upper bound, choose a subgraph $H\subseteq F$ attaining the
maximum in the definition of $\rho$, and put $h=v(H)$, $m=e(H)$.
Thus $\rho h=m+1$. Fix a set $S$ of size
\[
 k=\left\lceil\frac{[B+\log L]_+}{q}\right\rceil.
\]
Since $\rho\geq1$, we have $B\leq L$, and hence $k=o(n)$.
Also, in both cases $k=0$ and $k>0$,
\[
 n(1-p)^k\leq \frac{p^{-\rho}}L.
\]
Let $Z$ count the copies of $H$ among the vertices left undominated by
$S$. For any particular copy, its edges and its possible edges to $S$
are disjoint sets of random variables. Therefore
\[
 \bE [Z]
 \leq n^h(1-p)^{kh}p^m
 \leq \frac{p^{-\rho h+m}}{L^h}
 =\frac{1}{pL^h}.
\]
By Markov's inequality, whp $Z\leq1/p$. Add to $S$ one vertex from
each of these copies. The vertices left undominated by the enlarged set
form a subset of the original undominated set and contain none of its
$H$-copies. They are therefore $H$-free and, since $H\subseteq F$,
also $F$-free. This proves the upper bound in
\eqref{app:eq-general-F-bounds}.

For the lower bound, take a sufficiently large constant $C_F$ and put
\[
 k=\left\lfloor\frac{B-2\log L-C_F}{q}\right\rfloor.
\]
If the expression inside the floor is negative, the lower bound is immediate. Otherwise $k\geq0$, and
we show that whp no set of size $k$ is $F$-isolating.
For a fixed such set $S$, let $T=|U_S|$. Then
$T\sim\operatorname{Bin}(n-k,(1-p)^k)$ and, since
$k\leq L/p=o(n)$, its mean satisfies
\[
 \mu=(n-k)(1-p)^k
 \geq \frac{e^{C_F}}2p^{-\rho}L^2.
\]
Chernoff's inequality gives
\[
 \bP(T<\mu/2)\leq e^{-\mu/8}
 \leq\exp\left(-\frac{e^{C_F}}{16p}L^2\right).
\]
Conditional on $U_S$, the graph it induces is $G(T,p)$.
If $F$ has an edge, every subgraph $J\subseteq F$ with an edge satisfies
$v(J)\geq2$ and $\rho v(J)\geq e(J)+1$. Thus, whenever $T\geq\mu/2$,
\[
 T^{v(J)}p^{e(J)}
 \geq c_F p^{-\rho v(J)+e(J)}L^{2v(J)}
 \geq c_F p^{-1}L^4.
\]
Lemma~\ref{app:lem-janson-F} now bounds the conditional probability of
being $F$-free by $\exp(-c'_F L^4/p)$. If $F$ has no edges, the same
conditional probability is zero for all sufficiently large $n$, since
$T\geq\mu/2\to\infty$.

For $k\geq1$, the estimate used in Lemma~\ref{lem:expansion} gives
$\log\binom nk\leq L^2/p$; this also holds for $k=0$. A union bound
therefore bounds the probability of an $F$-isolating $k$-set by
\[
 \exp(L^2/p)\left[
   \exp\left(-\frac{e^{C_F}}{16p}L^2\right)
   +\exp(-c'_F L^4/p)\right]=o(1),
\]
where the second term may be omitted when $F$ has no edges.
Here it is enough to choose $e^{C_F}>32$.
Since adding vertices preserves $F$-isolation, this rules out all sets
of size at most $k$, proving the lower bound.

Finally, Lemma~\ref{log-1-p-approx} gives $1/p\leq1+1/q$.
The two bounds, together with $\iota(G(n,p),F)\geq0$, imply the
asserted additive estimate, including when $B\leq0$.
\end{proof}

\begin{cor}\label{app:cor-polynomial-F}
Let $F$ be a fixed nonempty graph and let $p=n^{-\alpha}$ for a fixed
$0<\alpha<1/\rho(F)$. Then whp
\[
 \iota(G(n,p),F)
 \sim \frac{(1-\alpha\rho(F))\log n}{p},
 \qquad
 \frac{\iota(G(n,p),F)}{\gamma(G(n,p))}
 \longrightarrow \frac{1-\alpha\rho(F)}{1-\alpha}.
\]
If $np^{\rho(F)}\to0$, then $p\,\iota(G(n,p),F)\to0$ in probability.
\end{cor}

\begin{proof}
Note $q\sim p$ and
$\log(n p^{\rho(F)})=(1-\alpha\rho(F))\log n$, whereas
$\log L=O(\log\log n)$. Apply Theorem~\ref{app:thm-general-F}, and
apply the same theorem with $F=K_1$ for the ratio. This gives the first assertion.
For the last one, take a density-maximizing subgraph $H$ as above
and let $Z$ count all its copies in $G(n,p)$. Choosing one vertex from
each copy gives an $F$-isolating set, and
\[
 p\,\bE [Z]\leq p n^h p^m=(np^{\rho(F)})^h\to0.
\]
Markov's inequality completes the proof.
\end{proof}

Some examples are listed below. In the middle two rows, $\ell\geq3$
and $r\geq2$. The final column is the endpoint of the range in which
Corollary~\ref{app:cor-polynomial-F} gives a positive coefficient of
$(\log n)/p$; it does not assert a limiting law at that endpoint.
\begin{center}
\begin{tabular}{c|c|c}
$F$ & $\rho(F)$ & $1/\rho(F)$\\ \hline
nonempty forest & $1$ & $1$\\[3pt]
$C_\ell$ & $1+1/\ell$ & $\ell/(\ell+1)$\\[3pt]
$K_r$ & $(r^2-r+2)/(2r)$ & $2r/(r^2-r+2)$\\[3pt]
$K_3$ & $4/3$ & $3/4$
\end{tabular}
\end{center}
For a cycle, every proper subgraph is a forest, so the whole cycle
attains the maximum. For a clique, a subgraph on $s$ vertices has at
most $\binom s2$ edges, and $(\binom s2+1)/s$ is nondecreasing for
integers $s\geq1$. In particular, for triangles and
$0<\alpha<3/4$ we obtain
\[
 \iota(G(n,n^{-\alpha}),K_3)
 \sim (1-4\alpha/3)n^\alpha\log n.
\]
Thus forests are precisely the fixed forbidden graphs for which
isolation and domination have the same first-order estimate throughout
the polynomial sparse range $0<\alpha<1$.

\begin{rmk}\label{app:rmk-finite-family}
The same estimate holds for a fixed finite nonempty family $\mc F$
of nonempty graphs, with
$\rho(\mc F):=\min_{F\in\mc F}\rho(F)$ in place of $\rho(F)$.
For the lower bound, use a member attaining this minimum.
For the upper bound, choose a density-maximizing subgraph $H_F$ for each
$F\in\mc F$, and after the initial choice of $S$, add one vertex
from each remaining copy of any $H_F$. Writing $h_F=v(H_F)$, its
expected number of copies is at most
\[
 p^{(\rho(F)-\rho(\mc F))h_F-1}L^{-h_F}
 \leq p^{-1}L^{-h_F}.
\]
There are finitely many members, so Markov's inequality gives
the same repair bound $1/p$ whp.
\end{rmk}

\begin{cor}\label{app:cor-extremal-F}
For every fixed nonempty graph $F$ and all sufficiently large integers
$d$, there is a graph $G$ with $\delta(G)\geq d$ such that
\[
 \frac{\iota(G,F)}{|V(G)|}
 \geq\frac{\log d-(\rho(F)+1)\log\log d-O_F(1)}{d}.
\]
\end{cor}

\begin{proof}
Take $n$ to be the least even integer at least $d\log d$, and put
$p=(d+d^{2/3})/n$. The Chernoff estimate in the proof of
Theorem~\ref{thm:extremal} again gives $\delta(G(n,p))\geq d$ whp.
Here $p\sim1/\log d$, $L=\log(np)=\log d+o(1)$, and
\[
 \log(n p^{\rho(F)})
 =\log d-(\rho(F)-1)\log\log d+o(1),
 \qquad nq=d\bigl(1+O(1/\log d)\bigr).
\]
Apply the lower bound in~\eqref{app:eq-general-F-bounds} and divide
by $n$. The change from $nq$ to $d$ costs only $O_F(1)$ in the
numerator, proving the claim.
\end{proof}

The regular construction in Section~\ref{sec:regular} uses a range in
which the coupling theorem of Gao, Isaev and McKay is sufficient.
We also mention the more recent theorem of Behague, Il'kovi\v{c} and
Montgomery~\cite{BIM}, which establishes the Kim--Vu sandwich conjecture
throughout the range $d\gg\log n$. Extending precise isolation estimates
to regular graphs would still require keeping track of the errors in
this coupling, since isolation itself is not monotone under adding edges.

\subsection{Sparse triangles are isolated one at a time}

For cyclic forbidden graphs, the estimate in Theorem~\ref{thm:sparse}
can fail by a polynomial factor. In a sufficiently sparse graph we can
say precisely what happens for triangles: almost every triangle requires
its own vertex in an isolating set.

\begin{prop}\label{prop:appendix-sparse-triangles}
Let $p=p(n)$ satisfy $np\to\infty$ and $n^4p^5\to0$.
Then, whp,
\[
 \iota(G(n,p),K_3)\sim\frac{(np)^3}{6}.
\]
In particular, if $p=n^{-\alpha}$ for a fixed $4/5<\alpha<1$,
then $\iota(G(n,p),K_3)\sim n^{3-3\alpha}/6$ whp.
\end{prop}

\begin{proof}
Let $Z$ count the triangles of $G\sim G(n,p)$. We have
$\bE Z=\binom n3p^3\sim(np)^3/6$. Two distinct triangle
indicators are independent unless the triangles share an edge, so
\[
 \operatorname{Var}(Z)=O(n^3p^3+n^4p^5)=o((\bE Z)^2).
\]
By Chebyshev's inequality, the number of triangles satisfies
$Z=(1+o(1))(np)^3/6$ whp.

For a triangle $T$, write $N_G[V(T)]$ for the union of the
closed neighborhoods of its three vertices. An isolating set for
triangles must meet $N_G[V(T)]$: otherwise all three vertices of $T$
remain undominated. Let $B$ count the pairs of distinct triangles
whose closed neighborhoods intersect. Pairs sharing a vertex
contribute $O(n^5p^6+n^4p^5)$ to $\bE B$. For two vertex-disjoint
triangles, an intersection of their closed neighborhoods requires
either an edge between the two triangles or a common neighbor
outside them. The union bound therefore gives
\[
 \bE B=O(n^5p^6+n^4p^5+n^6p^7+n^7p^8)
       =O(n^7p^8)=o(\bE Z),
\]
where we used $np\to\infty$ and $n^4p^5\to0$.
By Markov's inequality, $B=o(\bE Z)$ whp.

Delete at most one triangle for each of these $B$ pairs. The
remaining at least $Z-B$ triangles have pairwise disjoint closed
neighborhoods. Any triangle-isolating set consequently has at least
$Z-B$ vertices. Conversely, choosing one vertex from each triangle
always gives a triangle-isolating set. Hence $Z-B\leq\iota(G,K_3)\leq Z$,
which proves the assertion.
\end{proof}

\subsection{Arbitrary forbidden graphs at fixed edge density}

The comparison with domination remains particularly tight when
$p$ is fixed. This does not require a second-moment estimate for
$F$-isolating sets.

\begin{prop}\label{prop:appendix-dense-gap}
Let $F$ be a fixed graph with at least one vertex, and fix
$\eta\in(0,1/2)$. If $\eta\leq p=p(n)\leq1-\eta$, then
$0\leq\gamma(G(n,p))-\iota(G(n,p),F)\leq1$ whp.
\end{prop}

\begin{proof}
Write $q=-\log(1-p)$, $L=\log(np)$, and
$h_F(t)=\bP(G(t,p)\text{ is }F\text{-free})$.
The packing argument in Lemma~\ref{lem:large-sets-contain-f}
gives constants $c,C>0$, depending only on $F$ and $\eta$, such that
$h_F(t)\leq C e^{-ct^2}$ for every integer $t\geq0$.
If $F$ has no edges, this also holds because $h_F(t)=0$
for $t\geq|V(F)|$.
Consequently, for $x\geq0$,
\[
 \sum_{t\geq0}\frac{x^t}{t!}h_F(t)
 \leq C\sum_{t\geq0}\exp\bigl(t\log(x+1)-ct^2\bigr)
 \leq\exp\bigl(O_{F,\eta}(1+\log^2(x+1))\bigr),
\]
where the last inequality follows by completing the square and
bounding the resulting Gaussian sum uniformly in its center.

Let $Y_{F,k}$ count the $F$-isolating $k$-sets and put
$a=e^{-qk}$ and $x=(n-k)a/(1-a)$. As in
\eqref{eq:isolation-exact-probability}, exposing the edges incident
with a fixed $k$-set gives
\[
 \bE Y_{F,k}
 =\binom nk(1-a)^{n-k}
   \sum_{t=0}^{n-k}\binom{n-k}{t}
          \left(\frac a{1-a}\right)^t h_F(t).
\]
If $qk=L-2\log L+O(1)$, then $x=\Theta_\eta(L^2)$.
Using the preceding bound and the comparison with $f$ in
\eqref{eq:curves-accurate}, we obtain
\[
 \log\bE Y_{F,k}\leq f(k)+O_{F,\eta}((\log L)^2)=f(k)+o(L).
\]
Let $r=\lfloor\kappa_\gamma+1/2\rfloor$, with
$\kappa_\gamma$ as in Lemma~\ref{lem:centres}. Then
$r-1\leq\kappa_\gamma-1/2$ and
$r+1\geq\kappa_\gamma+1/2$. Since $f'(u)\geq L/2$
throughout the interval in question for all sufficiently large $n$,
we have $\bE Y_{F,r-1}\to0$ and $\bE X_{r+1}\to\infty$.
Markov's inequality rules out $F$-isolating sets of size $r-1$,
while Lemma~\ref{lem:domination-existence} gives a dominating set
of size $r+1$ whp. Thus
$r\leq\iota(G,F)\leq\gamma(G)\leq r+1$ whp.
\end{proof}

\subsection{The vertices left undominated}
\label{app:residual}

The potential $\Phi_q$ also describes the vertices left undominated by an
isolating set.  This gives a more precise explanation for the difference
between isolation and domination: a successful isolating set typically
leaves far fewer vertices undominated than an unconditioned set of the same
size, but still leaves substantially more than zero.

We first record a useful estimate for the potential.  The same function is
known as the \emph{deformed exponential}; in the notation of
Shapira and Tyomkyn~\cite{STpantograph}, our $\Phi_q(x)$ is
$f_{e^{-q}}(x)$.  In particular, it satisfies the pantograph equation
$\Phi_q'(x)=\Phi_q(e^{-q}x)$.  We give the estimates needed here directly,
including their uniformity when $q$ tends to zero.

For $z>0$, let $W(z)$ be the unique positive solution of
$W(z)e^{W(z)}=z$.  Thus $W$ denotes the positive real branch of the
Lambert $W$ function.  For $q,x>0$, define an integer-valued random variable
$T_{q,x}$ by
\[
 \bP(T_{q,x}=t)
 =\frac{x^t e^{-q\binom t2}}{t!\,\Phi_q(x)},\qquad t=0,1,2,\ldots.
\]

\begin{lem}\label{app:potential-saddle}
Fix $q_0>0$, and suppose that $0<q\leq q_0$ and $qx\to\infty$.
Put
\[
 v=W(qxe^{q/2}),\qquad m=v/q,\qquad A=q+1/m.
\]
Then, uniformly under these assumptions,
\begin{equation}\label{app:potential-lattice}
 \Phi_q(x)
 =(1+o(1))\frac{\exp\bigl((v^2+2v)/(2q)\bigr)}{\sqrt{2\pi m}}
       \sum_{j\in\mathbb Z}e^{-A(j-m)^2/2}.
\end{equation}
If, in addition, $q\to0$, this simplifies to
\begin{equation}\label{app:potential-sparse}
 \Phi_q(x)
 =(1+o(1))\frac{\exp\bigl((v^2+2v)/(2q)\bigr)}{\sqrt{1+v}},
 \qquad
 \sqrt q\,(T_{q,x}-m)\ \xrightarrow{\ d\ }\ N(0,1).
\end{equation}
For all $q,x>0$, the integer
$h=\lfloor W(qxe^q)/q\rfloor$ is a mode of $T_{q,x}$, and
\begin{equation}\label{app:potential-tail}
 \bP(|T_{q,x}-h|\geq s)
 \leq 2\exp\bigl(-q s(s-1)/2\bigr)
 \qquad(s=1,2,\ldots).
\end{equation}
\end{lem}

\begin{proof}
Write $w_t=x^t e^{-q\binom t2}/t!$.  The ratio
\[
 \frac{w_{t+1}}{w_t}=\frac{xe^{-qt}}{t+1}
\]
is at least one precisely when $t+1\leq W(qxe^q)/q$.
This proves the assertion about the mode.  It also shows that, for
$t\geq h$ and $s\geq1$,
\[
 \frac{w_{t+s}}{w_t}
 \leq \prod_{i=0}^{s-1}e^{-qi}
 =e^{-q s(s-1)/2}.
\]
Summing over $t\geq h$ bounds the upper tail by the right-hand side.
The same argument, using $w_{t-1}/w_t$ for $t\leq h$, bounds the lower
tail.  Their sum gives~\eqref{app:potential-tail}.

We next estimate the terms near their maximum.  For real $t>0$, put
\[
 F(t)=t(\log x+1-\log t)-\frac q2t(t-1).
\]
Its unique stationary point is $m$, since
$\log x+q/2=\log m+qm$.  At this point,
\[
 F(m)=\frac{v^2+2v}{2q},\qquad
 F''(m)=-q-1/m=-A,\qquad F'''(t)=1/t^2.
\]
The numbers $W(qxe^q)/q$ and $m$ differ by at most $1/2$:
differentiate the solution of $qz+\log z=b$ with respect to $b$,
and use $dz/db=(q+1/z)^{-1}\leq1/q$.
Consequently, $|h-m|\leq1$.

Consider integers $t$ with
$|t-m|\leq v^{1/4}/\sqrt q+2$.
Since $v\to\infty$, uniformly in this range we have $t/m\to1$.
Stirling's formula and Taylor's theorem therefore give
\begin{equation}\label{app:saddle-local}
 w_t=(1+o(1))\frac{e^{F(m)}}{\sqrt{2\pi m}}
                e^{-A(t-m)^2/2}.
\end{equation}
Indeed, the Taylor remainder is
$O(|t-m|^3/m^2)=O(\sqrt q\,v^{-5/4})+o(1)=o(1)$,
and the Stirling prefactor changes by a factor $1+o(1)$.
The tail bound~\eqref{app:potential-tail} shows that the terms outside
this range contribute $o(\Phi_q(x))$.
The same range contains a proportion $1-o(1)$ of the Gaussian lattice
sum in~\eqref{app:potential-lattice}: after multiplication by $\sqrt A$,
its endpoints tend to infinity, while the lattice spacing $\sqrt A$
is bounded.  Summing~\eqref{app:saddle-local} proves
\eqref{app:potential-lattice}.

If $q\to0$, then $A\to0$, and a Riemann sum gives
\[
 \sum_{j\in\mathbb Z}e^{-A(j-m)^2/2}
 =(1+o(1))\sqrt{2\pi/A}.
\]
Since $mA=v+1$, this proves the first assertion in
\eqref{app:potential-sparse}.  The same local estimate, summed over
intervals in $\sqrt A(t-m)$, proves convergence to the standard normal
law.  Finally, we have $A/q=1+1/v\to1$, so $\sqrt A$ can be replaced by
$\sqrt q$.
\end{proof}

The lattice sum in~\eqref{app:potential-lattice} is needed when $p$ stays
bounded away from zero.  In that case the fluctuations have bounded
size, and the position of $m$ relative to the integers can remain visible.
When $p\to0$, the fluctuation scale tends to infinity and these lattice
effects disappear.

\begin{prop}\label{app:conditional-residual}
Assume that $0<p\leq1-\eta$ and $np^2/\log^3n\to\infty$.
Let $S$ be a fixed $k$-set satisfying
$qk=L-2\log L+O(1)$, where $q=-\log(1-p)$ and $L=\log(np)$.
Put $\mu=(n-k)(1-p)^k$ and $v=W(q\mu e^{q/2})$.
Conditional on $S$ being isolating, the distribution of $|U_S|$ has
total variation distance $o(1)$ from the distribution of $T_{q,\mu}$.
In particular, under this conditioning,
\[
 |U_S|=(1+o_{\bP}(1))\frac{2\log L}{q}.
\]
If $p\to0$, then also
\[
 \sqrt q\left(|U_S|-\frac vq\right)
       \ \xrightarrow{\ d\ }\ N(0,1).
\]
\end{prop}

\begin{proof}
Set $N=n-k$ and $a=(1-p)^k$.  By conditioning on the undominated set,
the conditional probability of $|U_S|=t$ is proportional to
\[
 \binom Nt\left(\frac a{1-a}\right)^t e^{-q\binom t2}.
\]
For $t\leq20\log(q\mu)/q$, the estimates in the proof of
Lemma~\ref{lem:poissonization} give, uniformly,
\begin{align*}
 \binom Nt\left(\frac a{1-a}\right)^t e^{-q\binom t2}
 &=\frac{\mu^t}{t!}e^{-q\binom t2}
       (1-a)^{-t}\prod_{i=0}^{t-1}(1-i/N)\\
 &=(1+o(1))\frac{\mu^t}{t!}e^{-q\binom t2}.
\end{align*}
Here $t^2/N=o(1)$ and
$at=O(L^2\log L/(np^2))=o(1)$.
Lemma~\ref{lem:series} and the same Poissonization argument show that
the remaining terms have normalized mass $o(1)$ for both distributions.
Thus the two normalized distributions have total variation distance
$o(1)$.

Since $q\mu=\Theta(L^2)$, we have
$v=2\log L-\log(2\log L)+O(1)$.
The mode in Lemma~\ref{app:potential-saddle} differs from $v/q$ by
at most one.  Its tail estimate therefore yields
$T_{q,\mu}=v/q+O_{\bP}(q^{-1/2}+1)$, which proves the first
conclusion.  If $p\to0$, then $q\to0$, and the normal limit follows
from~\eqref{app:potential-sparse} and the total variation approximation.
\end{proof}

This statement concerns a fixed set conditioned on being isolating.
It also has a direct interpretation within a typical graph whenever
$M_k\to\infty$.  For every fixed $\varepsilon>0$, let $B_k$ count the
isolating $k$-sets for which
\[
 \left|\frac{q|U_S|}{2\log L}-1\right|>\varepsilon.
\]
The proposition gives $\bE B_k=o(M_k)$, whereas
Proposition~\ref{prop:second-moment} gives $Y_k/M_k\to1$ in probability.
Markov's inequality consequently gives $B_k/Y_k\to0$ in probability.
Thus all but a vanishing proportion of the isolating $k$-sets have the
residual size described above.  This argument does not assert the same
statement for every minimum isolating set.

\subsection{A more precise comparison of the two thresholds}
\label{app:gap-refined}

We retain the notation $\kappa_\gamma$ and $\kappa_\iota$ of
Lemma~\ref{lem:centres}.  The next estimate identifies a second term in
their difference and makes the rounding error explicit.

\begin{prop}\label{app:gap-second-order}
Suppose that $0<p\leq1-\eta$ and $\log(1/p)=o(\log n)$, and put
 $L=\log(np)$, $\ell=\log L$ and $q=-\log(1-p)$.
Then
\begin{equation}\label{app:gap-second-order-eq}
 \kappa_\gamma-\kappa_\iota
 =\frac{2}{qL^2}\bigl(\ell^2-\ell\log(2\ell)+O_\eta(\ell)\bigr).
\end{equation}
Moreover, whp, simultaneously for $(Z,\kappa)=(\gamma,\kappa_\gamma)$
and $(Z,\kappa)=(\iota,\kappa_\iota)$,
\begin{equation}\label{app:rounding-localization}
 \left\lceil\kappa-L^{-1}\right\rceil
 \leq Z(G(n,p))\leq
 \left\lceil\kappa+L^{-1}\right\rceil.
\end{equation}
\end{prop}

\begin{proof}
Write $\mu(u)=(n-u)e^{-qu}$ and
$H(u)=\log\Phi_q(\mu(u))$.  By Lemma~\ref{lem:centres}, both zeros
satisfy $qu=L-2\ell+O_\eta(1)$.
As $\log(q/p)=O_\eta(1)$, this gives
\[
 qu\log(en/u)
 =(L-2\ell+O_\eta(1))(L-\ell+O_\eta(1))
 =L^2+O_\eta(L\ell).
\]
At either zero we have $\mu(u)=u\log(en/u)$, with the additional
term $H(u)$ on the right at the isolation zero.
Lemma~\ref{lem:series} gives $qH(u)=O_\eta(\ell^2)$.
It follows that
\begin{equation}\label{app:mu-precise}
 q\mu(u)=L^2+O_\eta(L\ell)
\end{equation}
at both zeros and, by monotonicity of $\mu$, throughout the interval
between them.

For $R=\log(qx)\to\infty$, Lemma~\ref{app:potential-saddle} implies
\begin{equation}\label{app:potential-two-terms}
 \log\Phi_q(x)
 =\frac{R^2-2R\log R+O_\eta(R)}{2q},
 \qquad 0<q\leq-\log\eta.
\end{equation}
To see this, its saddle satisfies
$v+\log v=R+q/2$, and hence
$v=R-\log R+O_\eta(1)$.
The Gaussian lattice sum in~\eqref{app:potential-lattice} is bounded
above and below by positive constant multiples of $A^{-1/2}$,
uniformly for bounded $A$.
Thus its contribution together with $(2\pi m)^{-1/2}$ has logarithm
$O_\eta(\log(v+1)+1)$, which is absorbed by the error in
\eqref{app:potential-two-terms}.

At $u=\kappa_\iota$, equation~\eqref{app:mu-precise} gives
$R=2\ell+O_\eta(\ell/L)$.
Substitution in~\eqref{app:potential-two-terms} yields
\[
 H(\kappa_\iota)
 =\frac{2}{q}\bigl(\ell^2-\ell\log(2\ell)+O_\eta(\ell)\bigr).
\]
Also, uniformly between the zeros,
\[
 f'(u)=\log(n/u)+\bigl(q+(n-u)^{-1}\bigr)\mu(u)
      =L^2+O_\eta(L\ell).
\]
The mean value theorem, applied to
$f(\kappa_\gamma)-f(\kappa_\iota)=H(\kappa_\iota)$, proves
\eqref{app:gap-second-order-eq}.  The relative error
$O_\eta(\ell/L)$ from the derivative contributes only
$O_\eta(\ell^3/L)=o(\ell)$ inside the parentheses.

For the rounding assertion, both $f'$ and $(f+H)'$ are
$\Theta_\eta(L^2)$ throughout a fixed-size neighborhood of the
corresponding zero.  The assertion for $f'$ follows from its displayed
formula.  For the other derivative, Lemma~\ref{lem:series} gives
\[
 |H'(u)|
 =\bigl(q+(n-u)^{-1}\bigr)
      \frac{\mu(u)\Phi_q'(\mu(u))}{\Phi_q(\mu(u))}
 =O_\eta(\ell).
\]
The errors between these functions and the logarithms of the
expectations are $O(\log k)=o(L)$ by
\eqref{eq:curves-accurate}.  Consequently, the expectation tends to
zero at $k_-:=\lceil\kappa-L^{-1}\rceil-1$ and to infinity at
$k_+:=\lceil\kappa+L^{-1}\rceil$.
Markov's inequality excludes sets of size $k_-$, and
Lemma~\ref{lem:domination-existence} or
Proposition~\ref{prop:second-moment} supplies a set of size $k_+$.
This proves~\eqref{app:rounding-localization}.
\end{proof}

For example, let $p_*=(\log\log n/\log n)^2$ and suppose that
$p\sim cp_*$ for fixed $c>0$.
If $2/c$ is not an integer, then
\eqref{app:rounding-localization} implies, without any assumption on
the fractional parts of the two zeros, that whp
\[
 \gamma(G(n,p))-\iota(G(n,p))
 \in\{\lfloor2/c\rfloor,\lceil2/c\rceil\}.
\]
At the exceptional constants there is a small directional correction.
If $p=(2/j)p_*$ exactly, where $j$ is a fixed positive integer, then
\eqref{app:gap-second-order-eq} gives
\[
 \kappa_\gamma-\kappa_\iota
 =j\left(1-\frac{\log(2\log\log n)+O(1)}{\log\log n}\right).
\]
Its distance below $j$ is much larger than $L^{-1}$, and its distance
above $j-1$ stays bounded away from zero.  Thus whp the gap belongs to
$\{j-1,j\}$.  These conclusions locate the possible integers. The next subsection
determines the law of the isolation number whenever its expectation
threshold has a finite nonzero limit. A full distributional description
of the gap may additionally require joint information about the
dominating and isolating counts.

There is also a more accurate numerical version of the expectation
comparison. If $p\to0$, $\log(1/p)=o(\log n)$, and
$qk=L-2\log L+O(1)$, put $\mu=(n-k)e^{-qk}$ and
$v=W(q\mu e^{q/2})$. Stirling's formula, Lemma~\ref{lem:poissonization}
and~\eqref{app:potential-sparse} give
\[
 \bE X_k=(1+o(1))\frac{e^{f(k)}}{\sqrt{2\pi k}},
 \qquad
 M_k=(1+o(1))
 \frac{\exp\bigl(f(k)+(v^2+2v)/(2q)\bigr)}
      {\sqrt{2\pi k(1+v)}}.
\]
Indeed, $k^2/n=o(1)$ and $ne^{-2qk}=o(1)$ in this range, so
the finite-population and logarithmic errors are $o(1)$ in the exponent.
The prefactors matter when the expected count tends to a finite positive
number, even though they do not affect the first-order locations.

\subsection{A Poisson limit at the expectation threshold}
\label{app:poisson}

The last question in Section~\ref{sec:concluding-remarks} has an affirmative
answer throughout the range of Theorem~\ref{thm:variable-counting}.

\begin{thm}\label{thm:app-poisson}
Fix $\eta\in(0,1)$. Suppose that $0<p\leq1-\eta$,
$np^2/\log^3 n\to\infty$, $q=-\log(1-p)$ and
$qk=\log(np)-2\log\log(np)+O(1)$.
If $M_k\to\lambda\in(0,\infty)$, then $Y_k$ converges in distribution
to a Poisson random variable of mean $\lambda$. Consequently,
\[
 \bP\bigl(\iota(G(n,p))=k\bigr)\longrightarrow1-e^{-\lambda},
 \qquad
 \bP\bigl(\iota(G(n,p))=k+1\bigr)\longrightarrow e^{-\lambda}.
\]
\end{thm}

The useful observation is that conditioning a particular set $S$ to be
isolating has very little effect on the graph outside $S$. It plants an
independent set there, namely the vertices left undominated by $S$.
However, this planted set is much smaller than the independent sets
that already occur in the random graph. We first make this observation
precise.

For two probability distributions $\nu$ and $\nu'$ on the same finite
or countable space, write
$d_{\rm TV}(\nu,\nu')=\sup_A|\nu(A)-\nu'(A)|$ for their total variation
distance.

\begin{lem}\label{lem:app-planted}
Let $H\sim G(N,p)$, put $q=-\log(1-p)$, and let $t\leq N$ be a
nonnegative integer. Form $\widehat H$ by choosing a uniformly random
$t$-set $U$ and then sampling $G(N,p)$ conditioned on $U$ being independent.
If
\[
 \frac{t^2}{N}\longrightarrow0,
 \qquad
 \log\frac{N}{\max\{t,1\}}-\frac{qt}{2}\longrightarrow\infty,
\]
then $d_{\rm TV}(\mathcal L(H),\mathcal L(\widehat H))\to0$.
\end{lem}

\begin{proof}
The assertion is immediate for $t=0$. Let $Z_t$ be the number of
independent $t$-sets in $H$. The density of the law of $\widehat H$
relative to that of $H$ is $Z_t/\bE Z_t$. Hence Cauchy--Schwarz gives
\[
 d_{\rm TV}(\mathcal L(H),\mathcal L(\widehat H))
 \leq\frac12\sqrt{\frac{\bE Z_t^2}{(\bE Z_t)^2}-1}.
\]
If two uniformly chosen $t$-sets meet in $\ell$ vertices, conditioning
both to be independent saves precisely $\binom\ell2$ edge restrictions.
Thus
\[
 \frac{\bE Z_t^2}{(\bE Z_t)^2}
 =\sum_{\ell=0}^t
 \frac{\binom t\ell\binom{N-t}{t-\ell}}{\binom Nt}
 e^{q\binom\ell2}.
\]
For $\ell\geq1$, the summand is at most
\[
 \left(\frac{et^2}{N\ell}
       e^{q(\ell-1)/2}\right)^\ell.
\]
The logarithm of the expression in parentheses is convex as a function
of $\ell\in[1,t]$, so its maximum occurs at an endpoint. At $\ell=1$
the expression is $et^2/N=o(1)$; at $\ell=t$ it is
$(et/N)e^{q(t-1)/2}=o(1)$. The sum of all terms with $\ell\geq1$
is therefore $o(1)$, by comparison with a geometric series. The
second-moment ratio is at least one, so it tends to one, as required.
\end{proof}

We now prove Theorem~\ref{thm:app-poisson}. The count comparison 
below is included in detail as closeness in total variation alone
does not control the expectations of unbounded random variables.

\begin{proof}[Proof of Theorem~\ref{thm:app-poisson}]
Write $L=\log(np)$, $a=(1-p)^k$, $M=M_k$, and
$\rho=\bP(S\text{ isolates})$ for any fixed $k$-set $S$.
The estimates from Section~\ref{sec:concentration} give
\[
 k=\Theta_\eta(L/p),\qquad
 a=\Theta\left(\frac{L^2}{np}\right),\qquad
 ka=o(1),\qquad \frac{k^2}{n}=o(1).
\]
Put $b_n=\lceil40\log L/q\rceil$. Then $b_n=o(k)$,
$b_n^2/n=o(1)$, and $kb_n/n=o(1)$.

For $m\geq k$, let $\rho_m$ be the probability that a prescribed
$k$-set isolates $G(m,p)$, keeping $k$ and $p$ fixed. Deleting vertices
outside the prescribed set can only help isolation. Conversely, an
isolating set remains isolating after adding $h$ vertices if all the
new vertices have a neighbor in that set. These observations give the
exact inequalities
\begin{equation}\label{eq:app-deletion}
 \rho_n\leq\rho_{n-h}\leq(1-a)^{-h}\rho_n
 \qquad(0\leq h\leq n-k).
\end{equation}
In particular, uniformly for $h\leq k+b_n$,
$\rho_{n-h}=(1+o(1))\rho$.

Fix $S$, put $C=[n]\setminus S$, and let $W_S$ count the isolating
$k$-sets in $G[C]$. We first compare $W_S$ with $Y_k$ in an
unconditioned random graph. The expected number of isolating $k$-sets
meeting $S$ is at most $(k^2/n)M=o(1)$. For $T\subseteq C$, isolation
in $G$ implies isolation in $G[C]$, and the reverse implication holds
whenever $T$ dominates every vertex of $S$. The latter condition fails
with probability at most $ka$, independently of the event that $T$
isolates $G[C]$. By~\eqref{eq:app-deletion},
\[
 \bE W_S=\binom{n-k}{k}\rho_{n-k}
 \leq(1-a)^{-k}M=O(1).
\]
It follows that
\begin{equation}\label{eq:app-unconditional-count}
 \bP(Y_k\ne W_S)=o(1).
\end{equation}

Next condition on the event $E_S$ that $S$ is isolating, and let
$U=U_S\subseteq C$. Given $|U|=t$, the set $U$ is uniform among the
$t$-subsets of $C$. Given $U$, the graph on $C$ has distribution
$G(n-k,p)$ conditioned on $U$ being independent. The edges from $S$
to $U$ are absent, while for each vertex of $C\setminus U$ its edges
to $S$ are independent Bernoulli edges conditioned on at least one
being present. These choices are independent for different vertices,
and independent of the graph on $C$.

The exact probability formula~\eqref{eq:isolation-exact-probability},
together with the upper-tail bound in Lemma~\ref{lem:series} and
Lemma~\ref{lem:poissonization}, shows that
\begin{equation}\label{eq:app-residual-cutoff}
 \bP(|U|>b_n\mid E_S)=o(1).
\end{equation}
Indeed, with $x=(n-k)a/(1-a)$ we have $\log(qx)=2\log L+O_\eta(1)$,
and the terms with $t>b_n$ are bounded by the corresponding tail of
$\Phi_q(x)$.

Lemma~\ref{lem:app-planted} applies uniformly for $0\leq t\leq b_n$:
its first condition follows from $b_n^2/n=o(1)$, and its second follows
from
\[
 \log\frac{n-k}{b_n}-\frac{qb_n}{2}
 =L-O_\eta(\log L)\longrightarrow\infty.
\]
Averaging over $t$ and using~\eqref{eq:app-residual-cutoff}, we obtain
\begin{equation}\label{eq:app-exterior-tv}
 d_{\rm TV}\bigl(\mathcal L(G[C]\mid E_S),
                         \mathcal L(G[C])\bigr)=o(1).
\end{equation}

We must still check that, under this conditioning, the count $Y_k-1$
can be replaced by $W_S$. The subtracted one accounts for $S$ itself.
The estimates for positive intersection sizes in the proof of
Proposition~\ref{prop:second-moment} give
\begin{equation}\label{eq:app-overlapping-count}
 \sum_{\substack{T\ne S,\ |T|=k\\T\cap S\ne\emptyset}}
 \bP(E_T\mid E_S)=o(M)+o(1)=o(1).
\end{equation}
Here the near-diagonal pairs contribute $o(1)$ after conditioning,
and all other positive intersection sizes contribute $o(M)$.

It remains to compare the sets $T\subseteq C$ in the two counts.
Fix a possible residual set $U$ of size $t\leq b_n$.
For a given $T\subseteq C$, delete the vertices of $U\setminus T$.
Every edge still relevant to the isolation of $T$ is now an
unconditioned Bernoulli edge; the only forced absent edges that remain
have both endpoints inside $T$, and so are irrelevant. Therefore
\begin{equation}\label{eq:app-conditioned-candidate}
 \bP(T\text{ isolates }G[C]\mid E_S,U_S=U)
 \leq\rho_{n-k-t+|T\cap U|}
 \leq(1-a)^{-(k+b_n)}\rho=(1+o(1))\rho.
\end{equation}
In particular, the conditional expected value of $W_S$ is $O(1)$,
and the conditional expected number of its sets meeting $U$ is at most
 $(1+o(1))kt(n-k)^{-1}M=o(1)$.
For a set $T$ avoiding $U$ and a vertex $s\in S$, the conditional
probability that $s$ has no neighbor in $T$ is
\[
 \left(\frac{1-p-a}{1-a}\right)^k
 \leq(1-p)^k=a.
\]
This event depends only on edges from $S$ to $T$, and hence is
independent of the graph on $C$, conditional on $U$.
Consequently, the conditional expected number of sets counted by
$W_S$ that avoid $U$ but do not dominate $S$ is at most
$ka\,\bE[W_S\mid E_S,U_S=U]=o(1)$.
These estimates are uniform for $t\leq b_n$.
Together with~\eqref{eq:app-residual-cutoff}
and~\eqref{eq:app-overlapping-count}, they prove
 $\bP(Y_k-1\ne W_S\mid E_S)=o(1)$.
Combining this with~\eqref{eq:app-exterior-tv}
and~\eqref{eq:app-unconditional-count} gives
\begin{equation}\label{eq:app-size-bias}
 d_{\rm TV}\bigl(\mathcal L(Y_k-1\mid E_S),
                              \mathcal L(Y_k)\bigr)=o(1).
\end{equation}

For completeness, we finish the probabilistic argument directly.
All $k$-sets are interchangeable, so for every integer $j\geq0$,
\[
 \bP(Y_k-1=j\mid E_S)
 =\frac{j+1}{M}\bP(Y_k=j+1).
\]
This is the usual size-bias identity. Equation~\eqref{eq:app-size-bias}
therefore implies
\[
 \frac{j+1}{M}\bP(Y_k=j+1)-\bP(Y_k=j)\longrightarrow0.
\]
Since $\bE Y_k=M\to\lambda$, the distributions of $Y_k$ are tight.
Every subsequential limit has probabilities $\pi_j$ satisfying
$(j+1)\pi_{j+1}=\lambda\pi_j$. Their sum is one, so necessarily
$\pi_j=e^{-\lambda}\lambda^j/j!$. Thus the full sequence converges
to the Poisson distribution of mean $\lambda$.

Finally, the expectation comparison in
Theorem~\ref{thm:variable-counting} gives $M_{k-1}\to0$ and
$M_{k+1}\to\infty$. Hence $\iota(G(n,p))\in\{k,k+1\}$ whp,
and $\iota(G(n,p))\leq k$ precisely when $Y_k>0$. The conclusion follows.
\end{proof}

\printbibliography

\end{document}